\documentclass[]{amsart}

\usepackage{hyperref}
\usepackage{url}

\usepackage{booktabs}
\usepackage{algorithm}
\usepackage{amsmath}
\usepackage{amsfonts}
\usepackage{amssymb}
\usepackage{algorithmic}
\usepackage{array}
\usepackage{graphicx}
\usepackage{subcaption}
\usepackage[export]{adjustbox}
\usepackage{amsmath,amssymb,amsthm}
\usepackage{fancyvrb}
\usepackage{enumerate} 
\usepackage{cases} 
\usepackage{color} 
\usepackage{graphicx}
\usepackage{subcaption} 
\DeclareCaptionLabelFormat{blank}{}
\usepackage{bm}
\usepackage{xcolor}
\usepackage{mathrsfs}  

\usepackage{mathtools}
\mathtoolsset{showonlyrefs}

\usepackage{geometry}

\numberwithin{equation}{section} \numberwithin{figure}{section}

{\theoremstyle{remark} \newtheorem{remark}{Remark}}

{\theoremstyle{remark} \newtheorem{example}{Example}[section]}

\newcommand{\R}{\mathbb{R}}
\newcommand{\Rd}{{\mathbb{R}^d}}

\newcommand{\embh}[1]{\textbf{#1}}

\newcommand\restr[2]{{
  \left.\kern-\nulldelimiterspace 
  #1 
  \vphantom{\big|} 
  \right|_{#2} 
  }}

\begin{document} 
	
\title[Neural Very Weak Formulations]{Neural Very Weak Formulations enabling Hardware-Oriented deep PDE solvers.}
\author[G.~Acosta]{Gabriel Acosta}
\address[G.~Acosta]{Departamento de Matem\'atica, FCEyN, Universidad de Buenos Aires / IMAS, CONICET, Buenos Aires, Argentina}
\email{gacosta@dm.uba.ar}
\thanks{GA has been supported in part by PIP-2023, grant 11220220100246CO}

\author[F.~Bersetche]{Francisco ~Bersetche}
\address[F.~Bersetche]{Departamento de Matem\'atica, FCEyN, Universidad de Buenos Aires / IMAS, CONICET, Buenos Aires, Argentina}
\email{fbersetche@dm.uba.ar}
\thanks{FB has been supported in part  by PIP-2023, grant 11220220100246CO}

\begin{abstract}
We show, as a proof of concept, that least-squares very weak formulations of elliptic problems can be effectively discretized by neural networks possessing \emph{low regularity}, provided the test functions are drawn from appropriately smooth spaces. Apart from the immediate computational benefit of avoiding automatic differentiation, this approach, evaluated across various neural network spaces, demonstrates good performance even in challenging contexts, such as singular solutions and high dimensional settings.
Particular attention is paid to   trial functions based on \emph{step activations} and \emph{one bit} quantized linear functions, which are amenable to efficient \emph{hardware-oriented} implementations. 
 
\end{abstract}
 
\maketitle
\section{Introduction}
Originally introduced to handle nonsmooth solutions of elliptic equations, very weak formulations (VWF) trace back to the transposition method of Lions and Magenes \cite{LM} \cite[Section 2.52]{G}. In the finite element method (FEM), for instance, they have also proven  useful for the numerical treatment of problems with low regularity data. 
Beyond these well known applications  they provide a natural alternative to avoid automatic differentiation in neural networks (NN) based numerical methods for partial differential equations. Moreover,  they are suitable for the use of NN spaces with low regularity. Although these obvious particularities do not seem to have gone completely unnoticed, to the best of our knowledge, there are no studies exploring their potential applicability in detail.  Several factors might explain this fact. Differentiation based methods are more   natural and  there is an obvious reluctance to employ very weak formulations in general contexts, particularly if  smooth solutions are expected. Another disadvantage of these methods lies in the requirement of regular test functions -a computational challenge for some numerical methods such as the FEM  - and the fact that they, in its simplest form, 
offer error estimates in $L^2$ or even weaker norms. On the other hand, in recent years, deep learning has emerged as a powerful tool for approximating solutions to partial differential equations (PDEs) by minimizing loss functionals that enforce governing equations and boundary conditions \cite{E18,He20,liu2,liu1,PINN,DGM18,FNM}. Although not always competitive with classical methods in low-dimensional regimes, neural networks exhibit significant promise in high-dimensional contexts and on complex geometries \cite{weinan2022some,wojtowytsch2020can}. Also, formulations based on first-order systems \cite{DeepFOSLS,liu2,liu1,opschoor2024first} circumvent the need for second-order derivatives in the cost functional. However, these approaches impose still regularity requirements on the trial space, which precludes the use of discontinuous activation functions. Conversely, employing trial spaces with reduced regularity, such as those utilizing ReLU or step-function architectures, offers substantial computational advantages regarding energy consumption, processing speed, and memory efficiency. This is evidenced by the growing interest in quantized networks as a robust alternative for efficient computation. Ultimately, discontinuous and highly quantized neural representations are particularly well-suited for hardware-aware implementations, reducing computational complexity and aligning with the  constraints of emerging low-precision computing paradigm \cite{quant1, umuroglu2017finn}.

In this work, we propose a theoretical and computational framework grounded in a Neural Very Weak Formulation (NVWF), specifically tailored to this class of networks. Numerical results illustrate the potential of this approach; although not competitive on high precision architectures compared to more traditional NN methods, it enables a mathematically rigorous formulation for problems involving both smooth and singular data and solutions. More importantly, it demonstrates the viability of developing neural PDE solvers on binarized architectures.

This paper is organized as follows. In Section \ref{sec:vwf} we present the very weak residual least square  formulations that are the basis of our algorithms. Section \ref{sec:imp_details} introduce the implementation details and Section \ref{sec:numerical} is devoted to present numerical simulations.


\section{VWF  and NVWF}
\label{sec:vwf}
We begin this section by reviewing some elementary aspects  of very weak formulations (VWF).
Let $\Omega \subset \mathbb{R}^n$ be a bounded domain with boundary $\partial \Omega$. For simplicity, we consider the model problem: Find $u$ such that

\begin{equation}
\label{eq:poisson}
\left\lbrace
\begin{aligned}
    -\Delta u&=f,\, \text{in }\, \Omega\\ 
    u&=g,\,\text{on } \, \partial \Omega. 
\end{aligned}
\right.
\end{equation}

The regularity of solutions to \eqref{eq:poisson}  depends on both, the regularity of the data and the domain $\Omega$ \cite{G}. To simplify the exposition, we temporarily assume $\Omega$ smooth (or a convex polygon in $\R^2$). Consequently, if $f\in L^2(\Omega)$ and $g\in H^{3/2}(\partial \Omega)$, then $u\in H^2(\Omega)$ \cite{G}. 
Multiplying \eqref{eq:poisson} by a smooth test function $w\in H^2(\Omega)\cap H_0^1(\Omega)$ vanishing on $\partial \Omega$, and integrating by parts twice yields 
  \begin{equation}
      \label{eq:vweaknh}
      -\int_{\Omega}u\Delta w=-\int_{\Gamma}g\frac{\partial w}{\partial\eta}+\int_\Omega fw ,
  \end{equation}
where $\eta$ is the outward unit normal to $\Gamma$. For homogeneous Dirichlet conditions ($g=0$), this reduces to 
 \begin{equation}
      \label{eq:vweak}
      -\int_{\Omega}u\Delta w=\int_\Omega fw .
  \end{equation}
Notice that \eqref{eq:vweak} is well defined for any $u\in L^2(\Omega)$, motivating the VWF of \eqref{eq:poisson} which, for homogeneous boundaries, reads: 

(VWF) Find $u\in L^2(\Omega)$ satisfying \eqref{eq:vweak} for all $w\in H^2(\Omega)\cap H_0^1(\Omega)$. 

Well-posedness for this VWF follows from classical inf-sup arguments. Discrete approximations, like those given by the finite element method (FEM), similarly rely on discrete inf-sup conditions requiring tailor-made spaces that are often challenging to construct. The high regularity required for test functions further discourages a direct FEM approach.

Here, we employ neural network (NN) trial spaces. Because the natural solution space for the VWF is $L^2(\Omega)$,   NN spaces with low regularity suffice. For the test spaces, we use scaled translations of a smooth, compactly supported generating function $\varphi$ (in our case we take  $\varphi\in C^\infty_0(\R^n)$). Such translation invariant spaces, under appropriate assumptions on $\varphi$  are sufficiently rich \cite{MaSch}. In particular,  they can  approximate functions  defined in $\R^n$ up to a prescribed small saturation error,  even for a function $\varphi$ with a moderately large support. This crucial fact,  advantageously attenuates numerical errors during the stochastic evaluation of oscillatory integrals arising from integration by parts. On bounded domains, these properties require careful handling. For smooth $\Omega$, Sobolev extension theorems \cite{Evans} allow functions on $\Omega$ to be treated as restrictions of globally smooth functions. Building on this, we propose two distinct NN \embh{ansatzes}, (A1) and (A2), using test functions not strictly supported in $\Omega$. Both make use of a broadly supported generating test function within a residual least-squares formulation. Specifically, we define $\varphi$ as a fixed scaling of the standard mollifier 
\begin{equation}\label{eq:rho}
\rho(x) := 
\left\lbrace \begin{aligned}
    e^{-\frac{1}{1 - \|x\|^2}} \quad & \text{if } \|x\| \leq 1,\\
    0 \quad & \text{if } \|x\| > 1.\\
\end{aligned} \right.
\end{equation}
For a fixed parameter $h > 0$ (chosen to be relatively large), we set $\varphi (x)=\rho_h(x) := \frac{1}{h^n}\rho\left(\frac{x}{h}\right)$. The support of $\varphi$ is thus the closed ball $\bar{B}_h(0)$. Other   generating test functions  enjoying similar approximation properties  perform comparably. 

Hereafter, $v_\theta$ denotes a generic element of an NN space. Let $D_{\partial \Omega}$ and $D_{\partial \Omega}^S$ represent the distance and signed distance functions to $\partial \Omega$, respectively, with $D_{\partial \Omega}^S>0$ inside $\Omega$. We introduce surrogates $d_{\partial \Omega}$ and $d_{\partial \Omega}^S$, which are equivalent to their exact counterparts up to multiplicative constants. These surrogates can be chosen to be as smooth as $\Omega$ up to (or across) the boundary. 
For smooth $\Omega$, we may assume $\frac{\partial d_{\partial \Omega}}{\partial \eta}(x)=1$ on $\partial\Omega$, and  although $d_{\partial \Omega}:\bar\Omega\to \R$ cannot be smooth if $\Omega$ is non-smooth \cite{NikThom}, certain Lipschitz domains (e.g., piecewise smooth) admit a surrogate satisfying $\frac{\partial d_{\partial \Omega}}{\partial \eta}(x)=1$ almost everywhere on $\partial\Omega$. We exploit this fact in ansatz (A2) below.

For simplicity, we restrict our attention to the homogeneous case $g=0$ in our numerical experiments, the general case can be handled similarly. Here, Hardy's inequality \cite{Balinsky2015} guarantees that $\frac{u}{d_{\partial \Omega}}\in L^2(\Omega)$, even for solutions with low regularity  (e.g., for $u\in H_0^1(\Omega)$ only). In both (A1) and (A2), we seek approximate solutions to \eqref{eq:vweak} of the form $u_\theta=d_0 v_\theta$, where $d_0 \in \{d_{\partial \Omega}, d_{\partial \Omega}^S\}$. We see that shifting the unknown from $u_\theta$ to $v_\theta$ requires NN architectures with  $L^2(\Omega)$ approximation properties (a well-established feature for standard activation functions \cite{pinkus}). Remark \ref{rem:density} confirms this also holds for the one-bit quantized networks used herein.

\subsection{(A1): NVWF for smooth domains }
\label{subsec:A1}
 Our first ansatz operates similarly to the fictitious domain method (FDM) in the sense that assumes the solution defined outside  $\Omega$ with a strong imposition of the Dirichlet boundary condition (FDM typically employs a weak imposition via a Lagrange multiplier \cite{GPP} instead).   Then, calling $u_\theta=d_{\partial \Omega}^Sv_\theta$, we look for minimizers of 
 \begin{equation}
 \label{eq:fun1}
    \mathcal{L}(u_{\theta})=\int_{\Omega} \left[ \int_{B_h(x)} \tilde f(y) \cdot \rho_h(x-y) \, dy + \int_{B_h(x)} u_{\theta}(y) \cdot \Delta \rho_h(x-y) \, dy  \right]^2 \, dx
\end{equation}
where  $\tilde f$ stands for an  extension of $f$, preserving the smoothness of $f$ across $\partial \Omega$ (note that $B_h(x)\cap \Omega^c$ may be  nonempty).  We assume   $d^S_{\partial\Omega}$ is sufficiently smooth in the neighborhood of $\Omega$ given by $\tilde\Omega:=\{x\in \R^n: d^S_{\partial\Omega}(x)>-h \} $ (note that     $B_h(x)\subset \tilde \Omega$ for all $x\in \Omega$).
 Notice that equating the bracketed term in \eqref{eq:fun1} to zero implies, assuming a smooth enough $u_\theta$ and after integration by parts 
 $$\int_{B_h(x)} \left(\tilde f(y)  + \Delta u_{\theta}(y) \right)\cdot  \rho_h(x-y) \, dy=0$$
 which for a fixed $h$, and any $x\in \Omega$, establishes the orthogonality of the residual w.r.t. the space generated by translations of $\rho_h$.

\subsection{(A2): NVWF for more general domains}
If we assume the existence of a surrogate function $d_{\partial \Omega}$ for which   
$\frac{ \partial{d_{\partial \Omega}}}{\partial \eta}(x) =1$  a.e. along $\partial \Omega$ then, calling $u_\theta=d_{\partial \Omega}v_\theta$, we can look for minimizers of
\begin{equation}
\label{eq:fun2}
    \mathcal{J}(u_{\theta})=\int_{\Omega} \Bigg[ \int_{B_h(x) \cap \Omega} f(y) \cdot \rho_h(x-y) \, dy + \int_{B_h(x) \cap \Omega} u_{\theta}(y) \cdot \Delta \rho_h(x-y) \, dy 
\end{equation}
$$
+ \int_{B_h(x) \cap \partial \Omega} \rho_h(x-y) \cdot v_{\theta}(y) \, dS \Bigg]^2 \, dx
$$
Replacing the expression $u_\theta$ by $d_{\partial \Omega}v_\theta$, we get, at least formally and after integration by parts twice
$$
\int_{B_h(x) \cap \Omega} u_{\theta}(y) \cdot \Delta \rho_h(x-y) \, dy=\int_{B_h(x) \cap \Omega} \Delta u_{\theta}(y) \cdot  \rho_h(x-y) \, dy-\int_{B_h(x) \cap \partial \Omega} \rho_h(x-y) \cdot v_{\theta}(y) \, dS ,
$$
thanks to the fact  that $d_{\partial \Omega}$ vanishes identically and  $\frac{ \partial{d_{\partial \Omega}}}{\partial \eta}(x) =1$ a.e. $x$ on $\partial \Omega$, while  for $y\in \partial B_h(x)$, 
$\rho_h(x-y)=0=\frac{\partial\rho_h}{\partial \eta}(x-y)$. In particular, equating to zero the bracketed term in \eqref{eq:fun2}  implies
 $$\int_{B_h(x)\cap \Omega} \left(  f(y)  + \Delta u_{\theta}(y) \right)\cdot  \rho_h(x-y) \, dy=0,$$
 for fixed $h$ and any $x\in \Omega$, giving again an orthogonality of the residual w.r.t. translations of $\rho_h$.

\begin{remark}
Notice that in the case of smooth $\Omega$, both  (A1) and (A2) can be applied indistinctly. While our numerical tests show similar results in both cases, (A1) avoids the computation of the boundary term $\int_{B_h(x) \cap \partial \Omega} \rho_h(x-y) \cdot v_{\theta}(y) \, d S$.  

\end{remark}

\section{Implementation details} \label{sec:imp_details}

 As detailed above, our experiments employ the test space generated by translations of the function $\rho \in C_0^{\infty}(\Omega)$ defined in \eqref{eq:rho}. Observe that the Laplacian of $\rho$ is explicitly given by
\begin{equation}\label{eq:lap_rho}
\Delta \rho(x) = 
\left\lbrace \begin{aligned}
    \frac{ (8 - 2d)\|x\|^4 + 4(d - 1)\|x\|^2 - 2d}{(1 - \|x\|^2)^4}e^{-\frac{1}{1 - \|x\|^2}}  \quad & \text{if } \|x\| \leq 1,\\
    0 \quad & \text{if } \|x\| > 1,\\
\end{aligned} \right.
\end{equation}
which scales according to the relation $\Delta \rho_h(x) = \frac{1}{h^{d+2}}\Delta\rho\left(\frac{x}{h}\right)$. Given the oscillatory behavior of $\Delta \rho$, the numerical integration of terms involving this function requires careful treatment.

\subsection{Sampling strategy}\label{sec:sampling}

\begin{figure}[h] 
    \centering
    \includegraphics[width=0.4\textwidth]{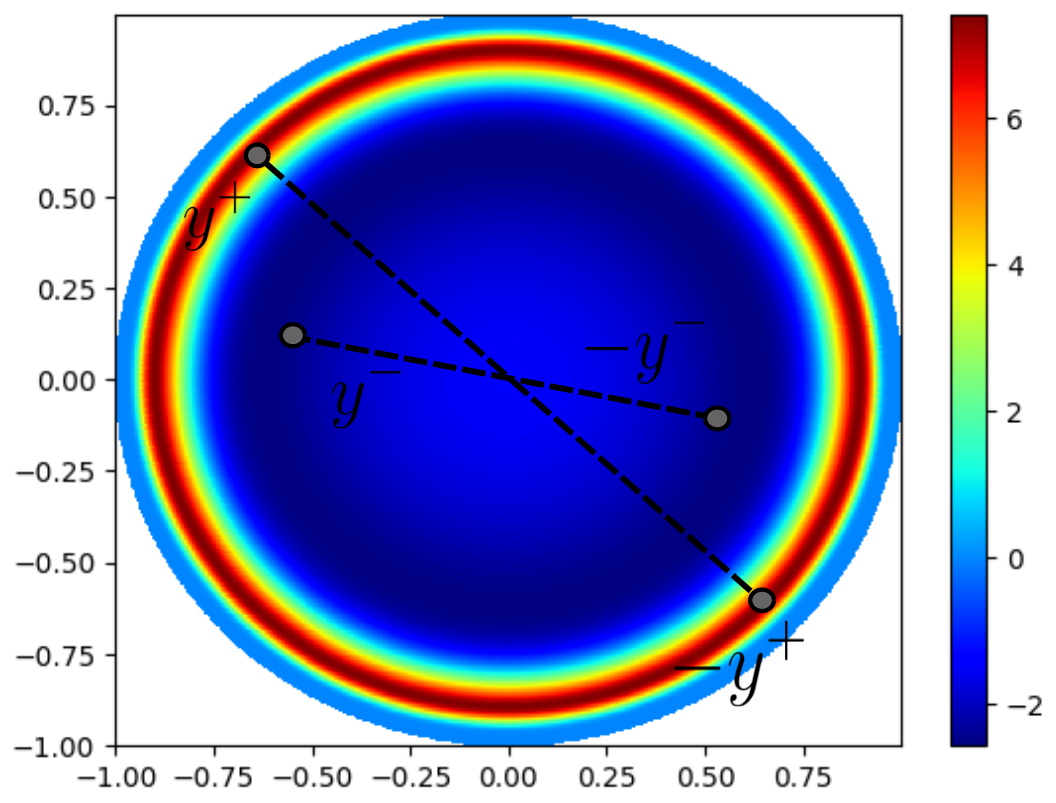}
\caption{Illustration of a representative sampling configuration based on the probability measures induced by $\Delta \rho$ for $N_y = 1$. This sampling strategy emulates a  centered difference operator.}
    \label{fig:sampling}
\end{figure}

The discretization of the functionals \eqref{eq:fun1} and \eqref{eq:fun2}  relies on approximating the involved integrals through Monte Carlo methods. While the exterior integral is computed using uniform sampling, the interior integrals are resolved using Monte Carlo integration augmented by importance sampling. Empirical results suggest that a balanced sampling strategy, specifically centered on the integral involving the Laplacian, improves the stability of the method. More precisely, setting $x = 0$ and $h = 1$ and letting $B = B_1(0)$ for brevity, we employ the following approximation:

\begin{equation}
\int_{B} u_{\theta}(y) \cdot \Delta \rho(-y) \, dy \approx \frac{ \int_{B} | \Delta \rho | }{4N_y} \Big( \sum_{y^+ \in Y^+}  u_{\theta}(y^+) + u_{\theta}(-y^+) - \sum_{y^- \in Y^-} u_{\theta}(y^-) + u_{\theta}(-y^-) \Big), 
\end{equation}
where $Y^+ \subset \{ \Delta \rho > 0 \}$ and $Y^- \subset \{ \Delta \rho < 0 \}$ are sets of points sampled according to the probability measures induced by the positive and negative parts of the Laplacian, $[\Delta \rho]^+$ and $[\Delta \rho]^-$, respectively. Here, $N_y$ denotes the cardinality of each set, being $N_y = \# Y^+ = \# Y^-$. By sampling an equal number of points in both regions and including their symmetric counterparts, we effectively emulate a   centered difference operator, as illustrated in Figure \ref{fig:sampling}.

\begin{remark}
    With a similar sampling strategy, one can approximate isolated terms such as $$\int_{B} u_{\theta}(y) \cdot\frac{\partial ^2\rho}{\partial^2  x_{i}} (-y).$$
     In particular, our method extends straightforwardly to more general second order problems, including those with variable coefficients.
\end{remark}

\subsection{Discretization of the loss functional $\mathcal{L}$}\label{sec:disc_func_L}
Following Section \ref{subsec:A1}, we call $\tilde{f}$ a smooth extension of the source term $f$. The trial solution is then constructed as $u_\theta = d_{\partial \Omega}^S v_\theta$, where $v_\theta$ represents a neural network parameterized by $\theta$, and $d_{\partial \Omega}^S$ stands for a surrogate signed distance function. Consequently, we define:  

\begin{equation}\label{eq:discrete_func_L}
    L(u_\theta) := \frac{|\Omega|}{N_x} \sum_{x \in X} \Bigg[   \frac{C_f}{N_f} \sum_{y \in Y^f_x } f(y) + \frac{ C_{\Delta} }{4N_y} \Big( \sum_{y^+ \in Y^+_x}  u_{\theta}(y^+) + u_{\theta}(-y^+) - \sum_{y^- \in Y^-_x} u_{\theta}(y^-) + u_{\theta}(-y^-) \Big) \Bigg]^2.
\end{equation}
As before, $X$ represents a set of $N_x$ points sampled uniformly over $\Omega$, and $Y^\Gamma$ denotes a set of $N_\Gamma$ points sampled uniformly over the boundary $\Gamma = \partial \Omega$. For each $x \in X$, the set $Y^f_x \subset B_h(x)$ is composed of $N_f$ points sampled according to the probability measure induced by $|\rho_h(x - \cdot)|$. Analogously, the sets $Y^+_x, Y^-_x \subset B_h(x)$ are generated following the strategy outlined in Section \ref{sec:sampling} by sampling from the measure induced by $|\Delta \rho_h(x - \cdot)|$, where $N_y = \#Y^+_x = \#Y^-_x$. The constants $C_f$ and $C_\Delta$ are consistent with the definitions provided in \eqref{eq:constants}.

\subsection{Discretization of the loss functional $\mathcal{J}$}\label{sec:disc_func_J}
The evaluation of the functional $\mathcal{J}$ requires computing integrals over the intersections $B_h(x) \cap \partial \Omega$ and $B_h(x) \cap \Omega$. To simplify  the numerical treatment of these terms, we perform the integration over  $\partial \Omega$ and $B_h(x)$ instead. This can be achieved  extending by zero the involved  integrands. Specifically, let $f_0: \mathbb{R}^n \to \mathbb{R}$ be the extension of $f$ defined as $f_0(x) = f(x)$ for $x \in \Omega$ and $f_0(x) = 0$ for $x \in \Omega^c$. We represent the trial solution as $u_\theta = d_0 v_\theta$, where $d_0$ is a surrogate distance function to the boundary $d_{\partial \Omega}$ extended by zero in $\Omega^c$ and $v_\theta$ is a neural network with parameters $\theta$. Since $\text{supp } \rho_h(x - \cdot) = B_h(x)$, the functional $\mathcal{J}$ in \eqref{eq:fun2} can be reformulated as:

\begin{equation}
    \mathcal{J}(u_{\theta})=\int_{\Omega} \Bigg[ \int_{B_h(x) } f_0(y) \cdot \rho_h(x-y) \, dy + \int_{B_h(x) } u_{\theta}(y) \cdot \Delta \rho_h(x-y) \, dy + \int_{ \partial \Omega} \rho_h(x-y) \cdot v_{\theta}(y) \, dS \Bigg]^2 \, dx.
\end{equation}

Building upon these observations, we define the discrete functional $J$, which approximates $\mathcal{J}$, as follows:

\begin{equation}\label{eq:discrete_func_J}
    J(u_\theta) := \frac{|\Omega|}{N_x} \sum_{x \in X} \Bigg[   \frac{C_f}{N_f} \sum_{y \in Y^f_x } f_0(y) + \frac{ C_{\Delta} }{4N_y} \Big( \sum_{y^+ \in Y^+_x}  u_{\theta}(y^+) + u_{\theta}(-y^+) - \sum_{y^- \in Y^-_x} u_{\theta}(y^-) + u_{\theta}(-y^-) \Big)
\end{equation}
$$ + \frac{|\Gamma|}{N_{\Gamma}} \sum_{y \in Y^\Gamma} \rho_h(x-y) \cdot v_\theta(y) \Bigg]^2.
$$
In this formulation, $X$ represents a set of $N_x$ points sampled uniformly over $\Omega$, while $Y^\Gamma$ denotes a set of $N_\Gamma$ points sampled uniformly over the boundary $\Gamma = \partial \Omega$. For every $x \in X$, the set $Y^f_x \subset B_h(x)$ contains $N_f$ points sampled according to the probability measure induced by $|\rho_h(x - \cdot)|$. Analogously, the sets $Y^+_x, Y^-_x \subset B_h(x)$ are generated following the strategy described in Section \ref{sec:sampling} by sampling from the measure induced by $|\Delta \rho_h(x - \cdot)|$, where $N_y = \#Y^+_x = \#Y^-_x$. The constants $C_f$ and $C_\Delta$ are given by:

\begin{equation} \label{eq:constants}
    C_f := \int_{B_h(0)} |\rho_h|, \quad C_\Delta := \int_{B_h(0)} |\Delta \rho_h|.
\end{equation}

\subsection{ Description of the method }
The proposed methodology focuses on approximating the minimizers of the functionals $J$ and $L$. Following the theoretical framework established in Section \ref{sec:vwf}, we utilize the functional $L$ specifically when the domain $\Omega$ is characterized by a smooth boundary. This approach enables us to bypass the numerical evaluation of the boundary flux term inherent in $J$. The implementation details for both variants of the proposed method are synthesized in Algorithms \ref{alg:flujo} and \ref{alg:dominio_ficticio}.


\begin{algorithm}[ht]
\caption{}
\label{alg:dominio_ficticio}
\begin{algorithmic}
    \STATE \textbf{Input:} $\Omega \subset \mathbb{R}^n$, $f \in L^{2}(\Omega)$, $d$ a surrogate signed distance function to $\partial \Omega$ (see Section \ref{sec:vwf}), $h \in \mathbb{R}$, $\rho_h \in C^{\infty}_0(\mathbb{R}^n)$ with $\text{supp}(\rho_h) = \overline{B_h(0)}$, and sample counts $N_x, N_f, N_y \in \mathbb{N}$.
    \STATE \textbf{Step $0$}: Define the smooth extension of $f$, denoted by $\tilde{f}$, as specified in Section \ref{sec:disc_func_L}. Initialize the neural network $v_\theta$ with parameters $\theta$ and define the trial solution $u_\theta = d v_\theta$.
    \WHILE{convergence criterion is not satisfied}
    \STATE \textbf{1:} Generate $X$ by sampling $N_x$ points uniformly from $\Omega$.
    \STATE \textbf{2:} For each $x \in X$:
    \begin{itemize}
    \item Generate $Y^f_{x}$ by sampling $N_f$ points in $B_h(x)$ following the probability measure induced by $|\rho_h(x - \cdot)|$.
    \item Generate $Y^+_{x}$ and $Y^-_{x}$ by sampling $N_y$ points each, following the balanced scheme described in Section \ref{sec:sampling}.
    \end{itemize}
    \STATE \textbf{3:} Compute the discrete functional $L$ according to \eqref{eq:discrete_func_L}.
    \STATE \textbf{4:} Update the network parameters $\theta$ using a gradient descent-based optimizer.
\ENDWHILE
\end{algorithmic}
\end{algorithm}


\begin{algorithm}[ht]
\caption{}
\label{alg:flujo}
\begin{algorithmic}
    \STATE \textbf{Input:} $\Omega \subset \mathbb{R}^n$, $f \in L^{2}(\Omega)$, $d$ a surrogate distance function to $\partial \Omega$ (see Section \ref{sec:vwf}), $h \in \mathbb{R}$, $\rho_h \in C^{\infty}_0(\mathbb{R}^n)$ with $\text{supp}(\rho_h) = \overline{B_h(0)}$, and sample counts $N_x, N_\Gamma, N_f, N_y \in \mathbb{N}$.
    \STATE \textbf{Step $0$}: Define the extensions of $f$ and $d$ outside $\Omega$ by zero, denoted by $f_0$ and $d_0$ respectively, as established in Section \ref{sec:disc_func_J}. Initialize the neural network $v_\theta$ with parameters $\theta$ and define the trial solution $u_\theta = d_0 v_\theta$.
    \WHILE{convergence criterion is not satisfied}
    \STATE \textbf{1:} Generate $X$ and $Y_\Gamma$ by sampling $N_x$ and $N_\Gamma$ points uniformly from $\Omega$ and $\partial \Omega$, respectively.
    \STATE \textbf{2:} For each $x \in X$:
    \begin{itemize}
    \item Generate $Y^f_{x}$ by sampling $N_f$ points in $B_h(x)$ following the probability measure induced by $|\rho_h(x - \cdot)|$.
    \item Generate $Y^+_{x}$ and $Y^-_{x}$ by sampling $N_y$ points each, following the balanced scheme described in Section \ref{sec:sampling}.
    \end{itemize}
    \STATE \textbf{3:} Compute the discrete functional $J$ according to \eqref{eq:discrete_func_J}.
    \STATE \textbf{4:} Update the network parameters $\theta$ using a gradient descent-based optimizer.
\ENDWHILE
\end{algorithmic}
\end{algorithm}

\subsection{Architectures with step activation functions}\label{subsec:step_arch}
We use the   notation 
$\mathcal{N}(L,P)$,
 to denote a generic neural network, with depth   $L\in \mathbb{N}$ (including the input and output layer)  and width vector  $P = (p_1, \dots, p_{L-2})$. Each $p_i\in \mathbb{N}$, for $1\le i\le L-2$, represent the width of layer $i$. We are interested in functions defined in $\R^d$, returning scalar values and therefore only the   width of hidden layers is declared.
 Furthermore, unless otherwise is stated, the simplified notation $\mathcal{N}(L,P)$ assumes that only Heaviside activations  $\sigma_s: \mathbb{R}^{p_l} \to \mathbb{R}^{p_{l}}$ are applied. In this context, it is well known that even architectures with two hidden layers, i.e.  $\mathcal{N}(4, (m,m))$,  can approximate functions in $L^p(\Rd)$, for $1\le p<\infty$   if we let $m\to \infty$ (see, for instance,  \cite{pinkus}). Regardless this theoretical result, for our implementation, we adopt the augmented Deep Heaviside Network (DHN) architecture proposed by \cite{kong2025expressivity}. This model extends the classical threshold network by integrating a small number of neurons with linear activations in each hidden layer alongside the standard Heaviside units. These linear components function as parallel information channels, allowing the network to preserve and propagate continuous input features throughout the deeper layers of the model.

The primary advantage of this hybrid design lies in its ability to overcome the representational bottlenecks typical of purely discrete activation functions. By combining the efficiency of Heaviside units with the expressivity of identity mappings, the architecture not only guarantees a richer feature propagation through deeper layers but also inherently mitigates the gradient flatlining problem during optimization, preserving active error signals across the continuous channels.

Formally, we denote the class of such networks as $\mathcal{N}(L, P, s)$, where $s$ is the number of augmented linear units per layer. The network is defined by the recursion $x_{l} = \sigma_s(W_l x_{l-1} + b_l)$, where the hybrid operator $\sigma_s: \mathbb{R}^{p_l + s} \to \mathbb{R}^{p_l + s}$ applies the Heaviside activation to the first $p_l$ units and the identity mapping to the remaining $s$ neurons.

Beyond its approximation properties, the localized discrete nature of the Heaviside units makes this architecture highly amenable to deployment on specialized digital accelerators (e.g., FPGAs). While the forward pass involves evaluating standard MAC operations for the linear bypass, the step activations drastically reduce the complexity of the nonlinear mapping overhead during the dense quadrature evaluations required by the very weak formulation.

The parameter optimization in our numerical experiments is performed via a surrogate gradient technique. Specifically, we employ the straight-through estimator (STE) with a sigmoid surrogate activation.
\begin{remark} 
\label{rem:density}
Let us notice that even extremely quantized versions of  $\mathcal{N}(4,(m,m))$ also retain density properties in $L^p$. Indeed, assuming weights undergo a 1-bit quantization constraining them to $\{-1, 1\}$, it is not difficult to prove  approximation properties if the bias vectors and the output layer weights belong to $\R$. Here we provide an elementary proof of this fact following ideas from \cite{pinkus}. We say that $K\subset \R^d$ is a $d$ -parallelotope with $(d-1)$- facets orthogonal to a certain base $B=\{v_1,\cdots,v_d\}$ of $\R^d$ if for each pair of parallel opposite facets there is an element in $B$ orthogonal to them. Recall that the number of facets of $K$ is $2d$. 

The key idea in the following argument is to show that there exists a \emph{fixed} basis of $\R^d$, $B=\{v_1,\cdots,v_d\}$ of quantized vectors such that 
the characteristic function of any  $d$ -parallelotope with $(d-1)$- facets orthogonal to $B$ belongs to $\mathcal{N}(4,(m,m))$. Taking into account that weights of the last layer are not quantized we can therefore obtain arbitrary linear combinations of such characteristic functions. It is not difficult to see that linear combinations of this type  approximate arbitrary functions in $L^p(\Rd)$ for large $m$. 

For our  quantized basis we choose (other choices are possible)   $v_1=(1,1,\cdots,1)$, and $v_i$ defined  as
$v_i=v_1-2e_i$ ($e_i$ is the $i-th$ canonical vector in $\R^d$) for each $2\le i\le d$. Note that 
all the components of each $v_i$ belong to the set $\{1,-1\}$. Calling $l_i(x)=v_i\cdot x+b_i$ and $m_i(x)=-v_i\cdot x+c_i$, where $\cdot$ stands for the scalar product and $b_i,c_i\in \R$, we see that \emph{any} $d$-parallelotope $K$ with  facets orthogonal $B$, can be expressed as
$$
K=\{x\in \Rd:l_i(x)>0 \quad m_i(x)>0, 1\le i\le d\},
$$
by choosing adequate biases $b_i,c_i\in \R$.
Therefore, the characteristic function $\kappa_K(x)$ of $K$ can be written, for instance, as
$$
\kappa_K(x)=\sigma\left(
(m-2d)\sigma(l_1(x))+\sum_{i=1}^d(\sigma(l_i(x))+\sigma(m_i(x)))-m+1/2\right)
$$
provided that $m\ge 2d$. Since the right hand side of the previous identity belongs to $\mathcal{N}(4,(m,m))$,  and weights of hidden layers belong to the set $\{-1,1\}$,  our claim follows. 
\end{remark}

\section{Computational experiments} \label{sec:numerical}

We evaluate the proposed methodology across several test cases. Numerical results demonstrate the robustness of the scheme under non-ideal regimes, successfully recovering the expected profiles even when the exact solutions exhibit low regularity. 

Optimization relies on the Adam algorithm across all computational configurations. To accelerate convergence, we employ a linear learning rate scheduler that decays monotonically to zero at the final epoch, with initial rates in the $[10^{-3}, 10^{-2}]$ range. We test our framework using distinct neural network trial spaces, denoted by $u_*$, where $* \in \{ \text{Tanh}, \text{ReLU}, \text{Step} \}$ specifies the activation function. The $\text{Step}$ variant designates a multi-layer feedforward network parameterized by discontinuous Heaviside operators. As formalized in Section~\ref{subsec:step_arch}, these architectures incorporate a fixed proportion of linear units within each hidden layer. Furthermore, to explore low-precision limits, some experiments employ a quantized architecture, $u_{\text{Qstep}}$, where hidden layer weights are constrained to 1-bit precision ($\{-1, 1\}$), while biases and output weights retain 32-bit floating-point precision. Network depth and width are specified for each respective test case. Finally, approximation accuracy is quantified via the error $e_* := u_* - u_{\text{GT}}$ in the $L^2$ norm, where $u_{\text{GT}}$ denotes the exact analytical solution.

The numerical simulations were implemented in PyTorch, and the source code is available on GitHub \footnote{\url{https://github.com/fbersetche/Neural-Very-Weak-Formulations-}}. The current software implementation of the architectures described in Section \ref{subsec:step_arch} and their quantized variants is executed on conventional hardware (CPUs/GPUs). Consequently, this work does not exploit low-level hardware-specific parallelism, and the present experiments are not intended to quantify execution efficiency. Nevertheless, our empirical findings demonstrate that these Heaviside-based models, in both single-precision and 1-bit quantized formats, achieve competitive approximation accuracy while maintaining training stability under gradient-based optimization. These results establish a numerical proof of concept for the proposed framework, validating its mathematical viability prior to deployment on dedicated low-precision hardware.

\begin{example}[1D] \label{sec:1D}

We evaluate the proposed method on the one-dimensional boundary value problem on $\Omega = [-1,1]$:
\begin{equation}\label{eq:1D}
\left\lbrace \begin{aligned}
	- u'' & = f \quad & \mbox{in } \Omega,\\
	u(-1) = u(1) & = 0.\\
\end{aligned} \right.
\end{equation}
To assess the scheme under varying regularity, we consider three source terms: $f_1=-\sin(\pi x)$, $f_2 = \delta$, and $f_3 = \delta'$, where $\delta$ and $\delta'$ denote the Dirac delta distribution centered at the origin and its distributional derivative, respectively. The corresponding exact solutions are $u_{\text{GT}} = \sin(\pi x)/\pi^2$, $u_{\text{GT}} = (1 - |x|)/2$, and $u_{\text{GT}} = (x+1)/2 - H(x)$, with $H: \R \to \R$ representing the Heaviside step function.

Numerical approximations are computed via Algorithm \ref{alg:dominio_ficticio}. The trial spaces for $u_{\text{ReLU}}$ are parameterized by NN featuring two hidden layers with $64$ neurons per layer. Similarly,  $u_{\text{Step}}$  uses the network topology $\mathcal{N}(4, (64,64), 6)$ defined in Section \ref{subsec:step_arch}. For low-precision regimes, the quantized architecture $u_{\text{QStep}}$ uses expanded configurations, $\mathcal{N}(4, (128, 128), 6)$, $\mathcal{N}(4, (256, 256), 6)$, and $\mathcal{N}(4, (512, 512), 6)$ for $f_1$, $f_2$, and $f_3$, respectively, under a 1-bit weight quantization, while biases and the output layer retain full 32-bit floating-point precision. Table \ref{table:1D} and Figure \ref{fig:1D} summarize the results and hyperparameters.

The method exhibits good performance across most test cases. However, for the lowest regularity solution ($f=\delta'$), the quantized architecture suffers from degraded convergence despite the increased network capacity, indicating that extreme quantization might induces spectral bias. This accuracy loss is mitigated by expanding the network dimensionality.
While the scheme demonstrates sensitivity to the choice of $h$, $N_x$, and $N_y$, optimizing these hyperparameter relations is left for future work. Another topic deserving attention it is the use of 
 post processing
 techniques that may improve in quality and accuracy the obtained quantized approximations. 

Next, we examine higher-dimensional settings where domain geometry introduces additional complexity.

\begin{table}[htbp]
    \centering
    \caption{Hyperparameter configurations ($h$, $N_x$, $N_y$) for the evaluated architectures and source terms in Example \ref{sec:1D} (1D experiments).}
    \label{table:1D}
    \begin{tabular}{l ccc ccc ccc}
        \toprule
        & \multicolumn{3}{c}{$f = -\sin(\pi x)$} & \multicolumn{3}{c}{$f = \delta$} & \multicolumn{3}{c}{$f = \delta'$} \\
        \cmidrule(lr){2-4} \cmidrule(lr){5-7} \cmidrule(lr){8-10}
        Architecture & $h$ & $N_x$ & $N_y$ & $h$ & $N_x$ & $N_y$ & $h$ & $N_x$ & $N_y$ \\
        \midrule
        $u_{\text{ReLU}}$   &  0.2   &  80   &  40  &  0.2   &   80  &  40   &  0.2   &  80   &  40   \\
        $u_{\text{Step}}$   &   0.2  &  80   &  40  &  0.2   &   80  &  40   &  0.2   &  80   &  80   \\
        $u_{\text{QStep}}$  &  0.2   &  80   &  40  &  0.2   &   40  &  100  &  0.2   &  100  &  160  \\
        \bottomrule
    \end{tabular}
\end{table}

\begin{figure}[htbp]
    \centering

    \begin{subfigure}{0.3225\textwidth}
        \centering
        \includegraphics[width=\linewidth]{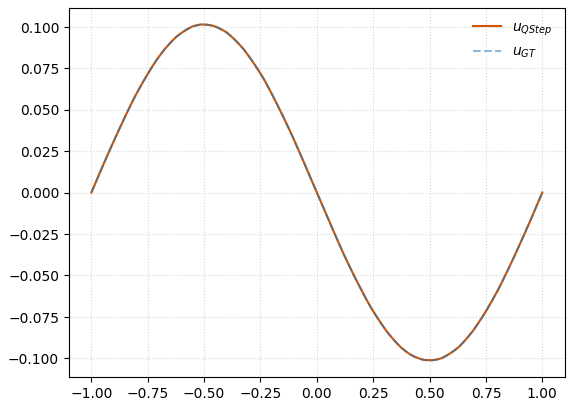}

    \end{subfigure}\hspace{0.5em}
    \begin{subfigure}{0.31\textwidth}
        \centering
        \includegraphics[width=\linewidth]{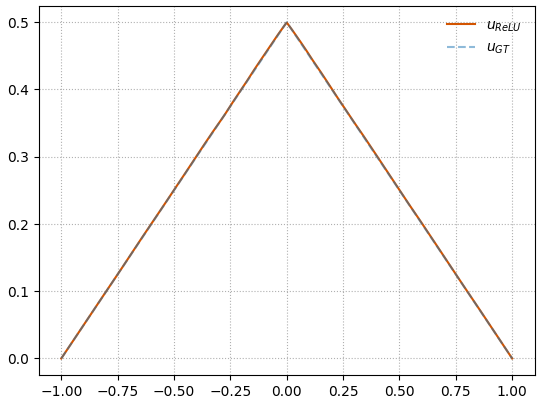}

    \end{subfigure}\hspace{0.5em}%
    \begin{subfigure}{0.315\textwidth}
        \centering
        \includegraphics[width=\linewidth]{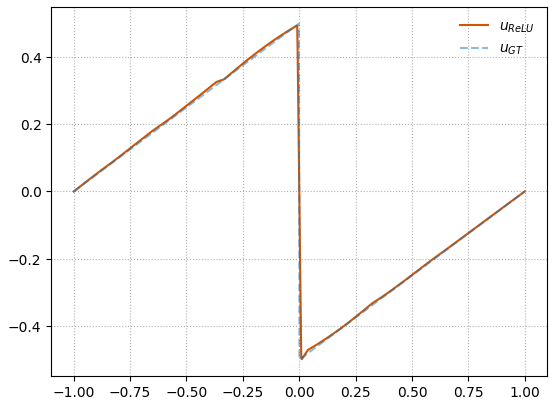}

    \end{subfigure}

    \vspace{0.5em}

    
    \begin{subfigure}{0.3225\textwidth}
        \centering
        \includegraphics[width=\linewidth]{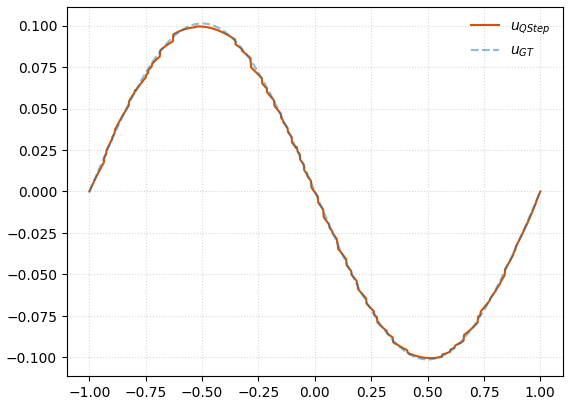}

    \end{subfigure}\hspace{0.5em}
    \begin{subfigure}{0.31\textwidth}
        \centering
        \includegraphics[width=\linewidth]{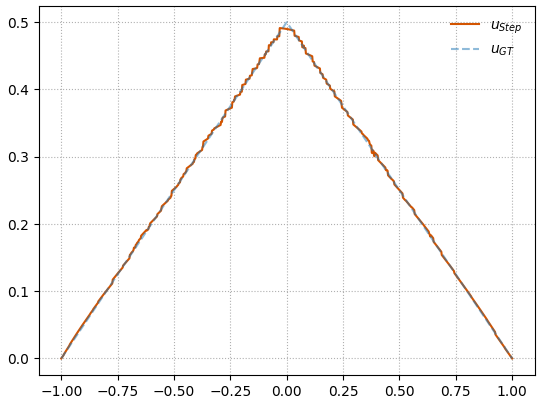}

    \end{subfigure}\hspace{0.5em}
    \begin{subfigure}{0.315\textwidth}
        \centering
        \includegraphics[width=\linewidth]{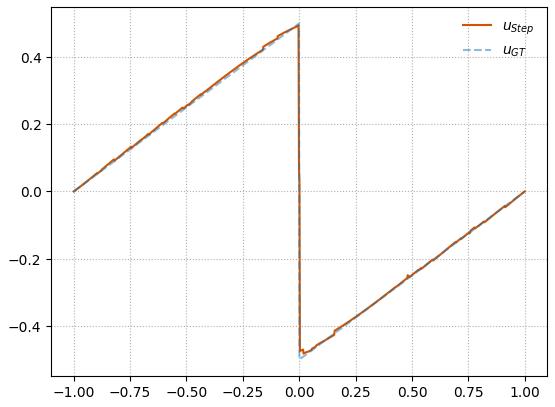}

    \end{subfigure} 

    \vspace{0.5em} 


    \begin{subfigure}{0.3225\textwidth}
        \centering
        \includegraphics[width=\linewidth]{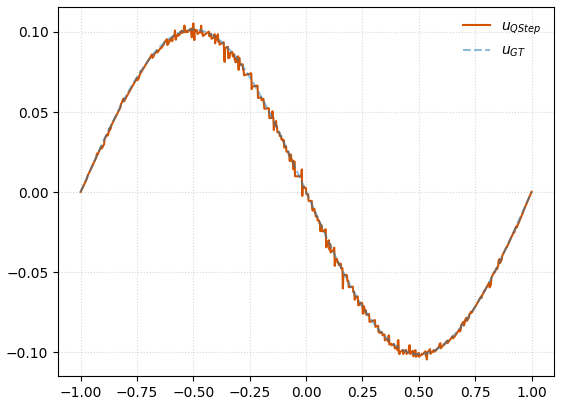}

    \end{subfigure}\hspace{0.5em}
    \begin{subfigure}{0.31\textwidth}
        \centering
        \includegraphics[width=\linewidth]{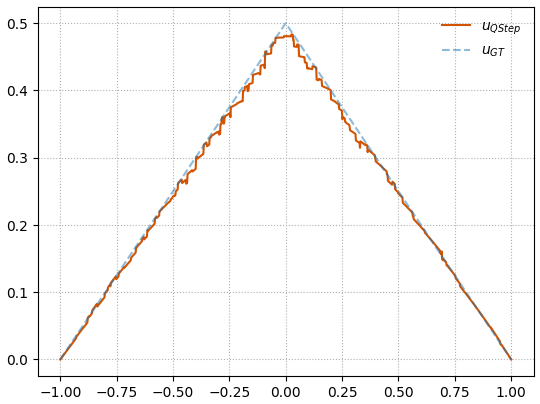}

    \end{subfigure}\hspace{0.5em}%
    \begin{subfigure}{0.315\textwidth}
        \centering
        \includegraphics[width=\linewidth]{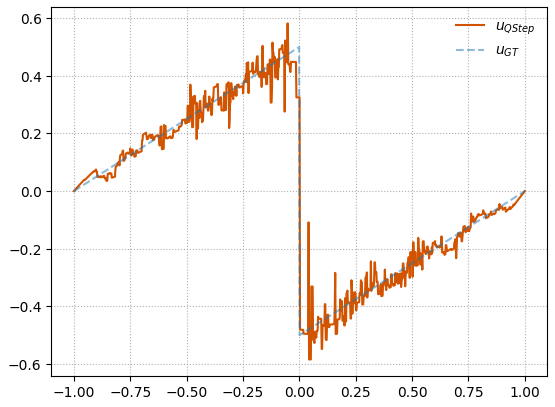}

    \end{subfigure} 

    \caption{  Numerical approximations (solid orange) and exact solutions $u_{\text{GT}}$ (dashed light blue) for the 1D test cases in Example \ref{sec:1D}. Columns (left to right) correspond to source terms $f = -\sin(\pi x)$, $f = \delta$, and $f = \delta'$. Rows (top to bottom) correspond to architectures $u_{\text{ReLU}}$, $u_{\text{Step}}$, and $u_{\text{QStep}}$.  }
    \label{fig:1D}
\end{figure}

\end{example}

\begin{example}[Smooth Solution - Smooth Domain]\label{sec:suave2}

We consider the unit ball $\Omega = B_1(0) = \{ x \in \mathbb{R}^2 : \|x\| < 1 \}$ and seek a function $u: \Omega \to \mathbb{R}$ such that

\begin{equation}\label{eq:suave2}
\left\lbrace \begin{aligned}
	-\Delta u & = -9\|x\| \quad & \mbox{in } \Omega,\\
	u & = 0 \quad & \mbox{on } \partial\Omega.\\
\end{aligned} \right.
\end{equation}
The solution to \eqref{eq:suave2} is
\[
u_{\text{GT}} = \|x\|^3 - 1.
\]

To approximate the solution to this problem, we implement the algorithmic procedure outlined in Algorithm \ref{alg:dominio_ficticio}. The functional spaces for all baseline trial functions are parameterized via NN consisting of two hidden layers with an internal dimensionality of $128$ neurons per layer. Specifically, $u_{\text{Step}}$ incorporates a network  topology $\mathcal{N}(4, (128,128), 4)$ (see   Section \ref{subsec:step_arch}). To explore low-precision operational regimes, this example further incorporates a  1-bit weight quantization network $u_{\text{QStep}}$ $\mathcal{N}(4, (256, 256), 4)$ configuration. 

Numerical results and  the selected hyperparameter configurations, are summarized in Table \ref{table:suave2} and visualized in Figure \ref{fig:suave2}.

\begin{table}[htbp]
\centering
\caption{Hyperparameters and resulting relative $L^2$ errors for Algorithm \ref{alg:dominio_ficticio} applied to Example \ref{sec:suave2} (Smooth Solution - Smooth Domain).}
\label{table:suave2}
\begin{tabular}{l ccc r}
\toprule
Architecture & $h$ & $N_x$ & $N_y$ & $\|e_*\|_2 / \|u_{\text{GT}}\|_2$ \\ 
\midrule
$u_{\text{Tanh}}$ & 0.2 & 20 & 100 & 0.0038 \\ 
$u_{\text{ReLU}}$ & 0.2 & 20 & 100 & 0.0022 \\ 
$u_{\text{Step}}$  & 0.2 & 100 & 20 & 0.013 \\ 
$u_{\text{QStep}}$  & 0.2 & 40 & 100 & 0.043 \\ 
\bottomrule
\end{tabular}
\end{table}

\begin{figure}[htbp]
    \centering
    
    \begin{subfigure}{0.22\textwidth}
        \centering
        \includegraphics[width=\linewidth]{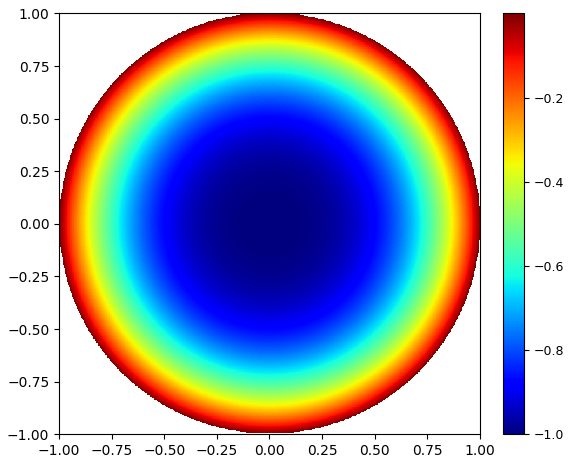}
        \caption{$u_{\text{GT}}$}
    \end{subfigure}\hfill
    \begin{subfigure}{0.22\textwidth}
        \centering
        \includegraphics[width=\linewidth]{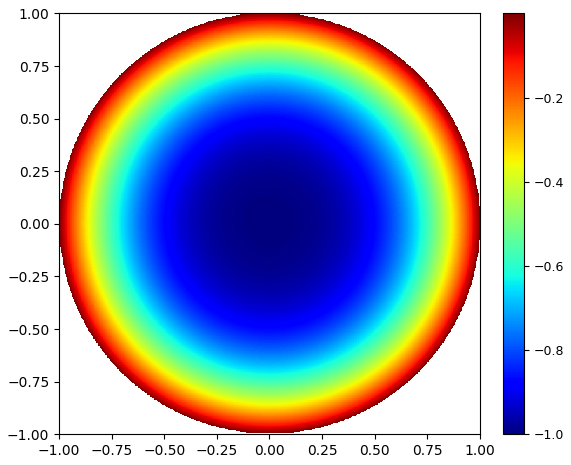}
        \caption{$u_{\text{Tanh}}$}
    \end{subfigure}\hfill
    \begin{subfigure}{0.24\textwidth}
        \centering
        \includegraphics[width=\linewidth]{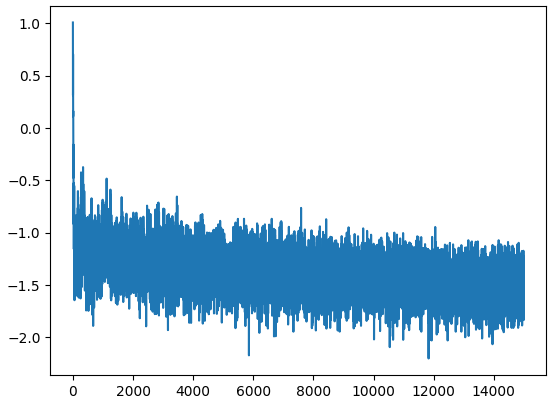}
        \caption{$\log(loss)$}
    \end{subfigure}\hfill
    \begin{subfigure}{0.24\textwidth}
        \centering
        \includegraphics[width=\linewidth]{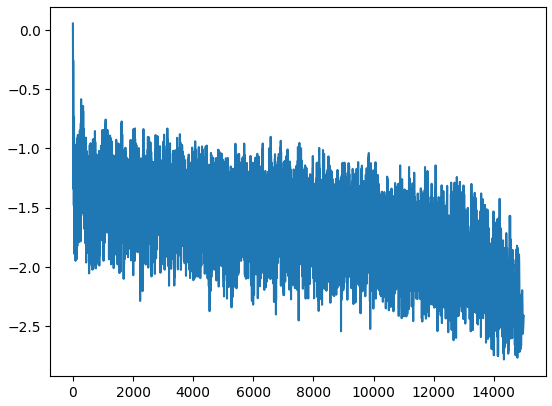}
        \caption{$\log( \|e_{\text{Tanh}}\|_{2} / \|u_{\text{GT}}\|_{2} )$}
    \end{subfigure}

    \vspace{1em} 

    \begin{subfigure}{0.22\textwidth}
        \centering
        \includegraphics[width=\linewidth]{suave2_u_GT.png}
        \caption{$u_{\text{GT}}$}
    \end{subfigure}\hfill
    \begin{subfigure}{0.22\textwidth}
        \centering
        \includegraphics[width=\linewidth]{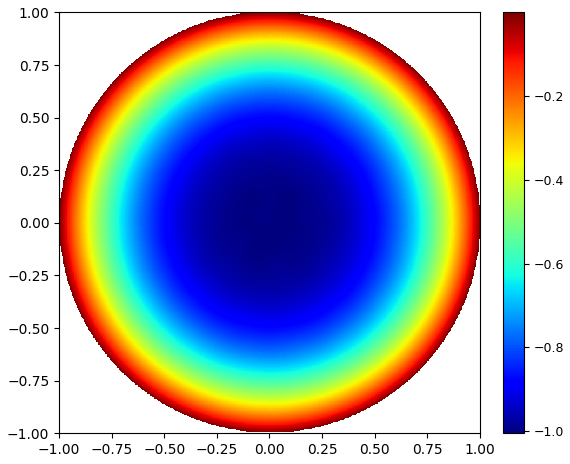}
        \caption{$u_{\text{ReLU}}$}
    \end{subfigure}\hfill
    \begin{subfigure}{0.24\textwidth}
        \centering
        \includegraphics[width=\linewidth]{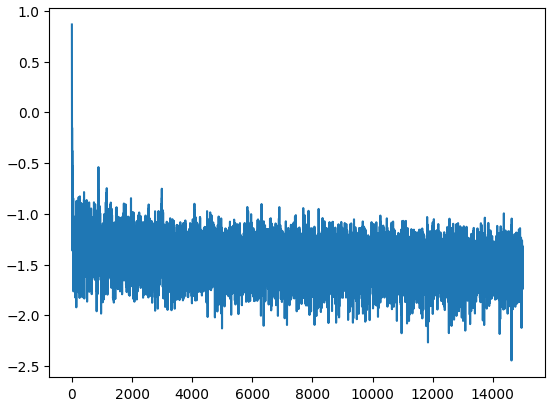}
        \caption{$\log(loss)$}
    \end{subfigure}\hfill
    \begin{subfigure}{0.24\textwidth}
        \centering
        \includegraphics[width=\linewidth]{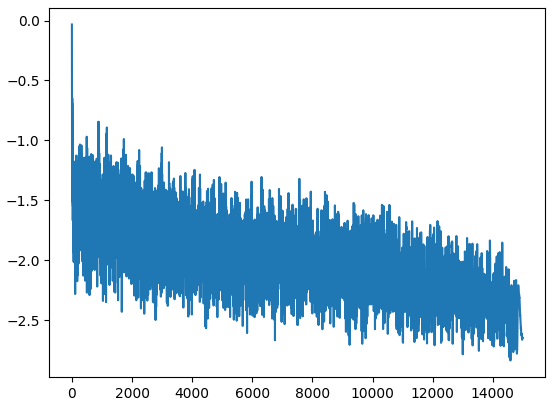}
        \caption{$\log(\|e_{\text{ReLU}}\|_{2} / \|u_{\text{GT}}\|_{2})$}
    \end{subfigure}

    \vspace{1em}

    \begin{subfigure}{0.22\textwidth}
        \centering
        \includegraphics[width=\linewidth]{suave2_u_GT.png}
        \caption{$u_{\text{GT}}$}
    \end{subfigure}\hfill
    \begin{subfigure}{0.22\textwidth}
        \centering
        \includegraphics[width=\linewidth]{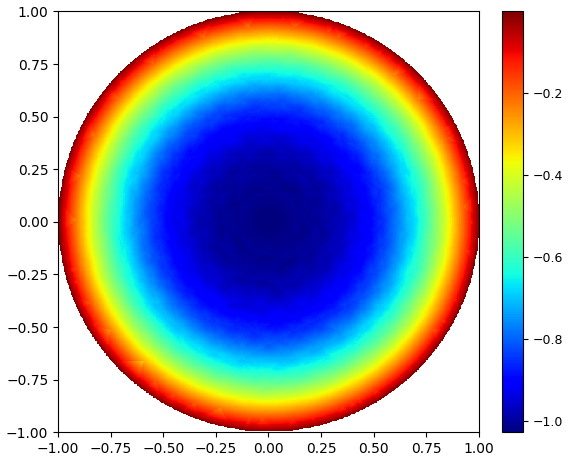}
        \caption{$u_{\text{Step}}$}
    \end{subfigure}\hfill
    \begin{subfigure}{0.24\textwidth}
        \centering
        \includegraphics[width=\linewidth]{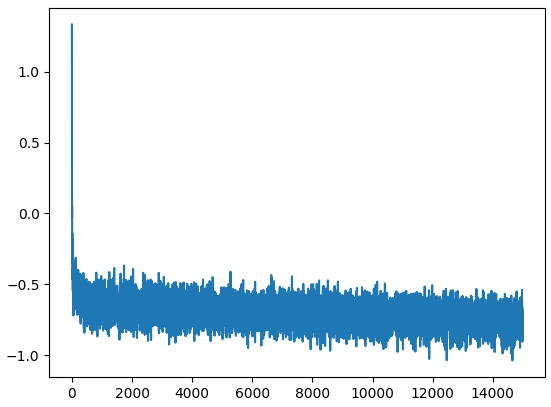}
        \caption{$\log(loss)$}
    \end{subfigure}\hfill
    \begin{subfigure}{0.24\textwidth}
        \centering
        \includegraphics[width=\linewidth]{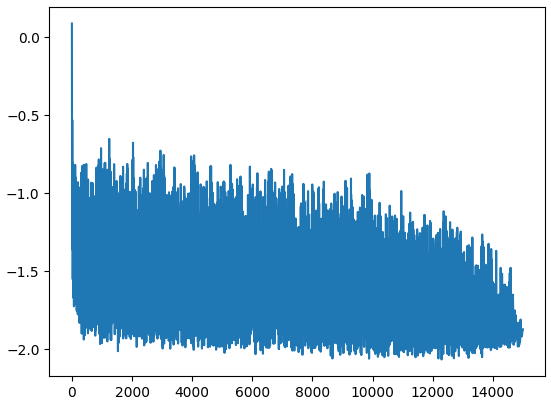}
        \caption{$\log(\|e_{\text{Step}}\|_{2} / \|u_{\text{GT}}\|_{2})$}
    \end{subfigure}

        \begin{subfigure}{0.22\textwidth}
        \centering
        \includegraphics[width=\linewidth]{suave2_u_GT.png}
        \caption{$u_{\text{GT}}$}
    \end{subfigure}\hfill
    \begin{subfigure}{0.22\textwidth}
        \centering
        \includegraphics[width=\linewidth]{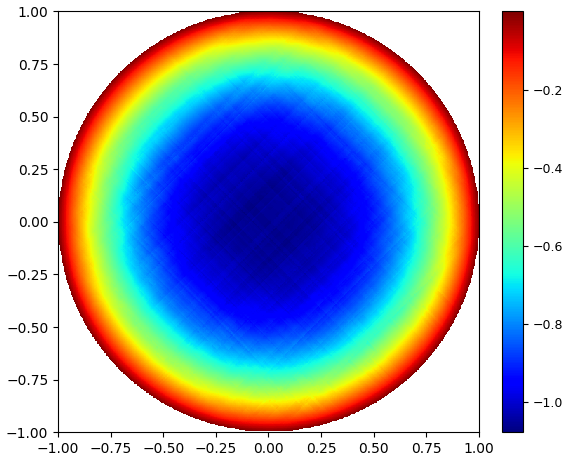}
        \caption{$u_{\text{QStep}}$}
    \end{subfigure}\hfill
    \begin{subfigure}{0.24\textwidth}
        \centering
        \includegraphics[width=\linewidth]{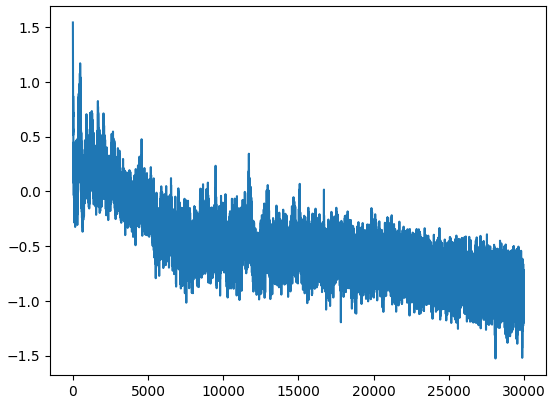}
        \caption{$\log(loss)$}
    \end{subfigure}\hfill
    \begin{subfigure}{0.24\textwidth}
        \centering
        \includegraphics[width=\linewidth]{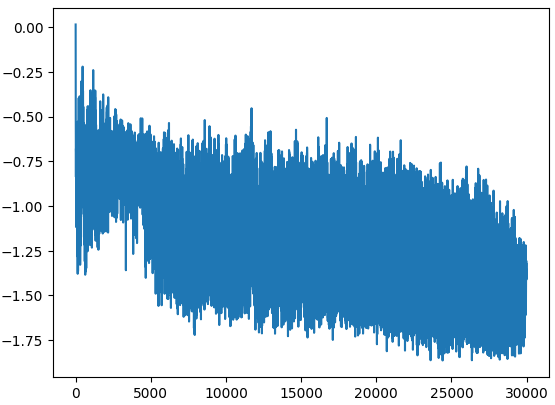}
        \caption{$\log(\|e_{\text{QStep}}\|_{2} / \|u_{\text{GT}}\|_{2})$}
    \end{subfigure}

    \vspace{1em}

    \caption{Performance of four neural network architectures for Example \ref{sec:suave2} using Algorithm \ref{alg:dominio_ficticio}. Rows correspond to architectures. Columns (left to right): exact solution, numerical approximation, loss history (log scale), and relative $L^2$ error history (log scale).}
    \label{fig:suave2}
\end{figure}

\end{example}

\begin{example}[Smooth Solution - Non-smooth Domain]\label{sec:suave}
We consider the following 2-dimensional problem. Let $\Omega = [-1,1]^2,$ and $k \in \mathbb{N}$. We seek $u \colon \Omega \to \R$ such that
\begin{equation}\label{eq:example1}
\left\lbrace \begin{aligned}
	-\Delta u & = \sin( \pi x_1) \left( \pi^2 (1-x_2^2) + 2 \right) \quad & \mbox{in } \Omega,\\
	u & = 0 \quad & \mbox{on } \partial\Omega.\\
\end{aligned} \right.
\end{equation}
The solution to \eqref{eq:example1} is
\[
u_{\text{GT}} = \sin( \pi x_1)(1-x^2_2).
\]

Due to the non-smooth character of the boundary $\partial \Omega$, this example is treated using the framework described in Algorithm \ref{alg:flujo}. In alignment with our previous setup, the approximation spaces for all trial functions are parameterized via multilayer neural networks featuring two hidden layers of $128$ neurons each. For $u_{\text{Step}}$ integrates a neural network $\mathcal{N}(4, (128,128), 16)$ (see Section \ref{subsec:step_arch}). This example is also tested evaluating a quantized variant of the step-activation architecture, denoted as $u_{\text{QStep}}$, using a network configuration $\mathcal{N}(4, (2048,2048), 4)$. It is worth noting that this network topology is substantially larger than the one deployed in our previous experiments. This architectural expansion compensates 
 the representation capacity degradation induced by  1-bit quantization when   dealing with  more complex  solutions.

 Table \ref{table:suave} and  Figure \ref{fig:suave} summarize the obtained numerical results and hyperparameter configurations.

\begin{table}[htbp]
\centering
\caption{Hyperparameters and resulting relative errors for Algorithm \ref{alg:flujo} applied to Example \ref{sec:suave} (Smooth Solution - Non-smooth Domain).}
\label{table:suave}
\begin{tabular}{l cccc r}
\toprule
Architecture & $h$ & $N_x$ & $N_y$ & $N_b$ & $\|e_*\|_2 / \|u_{\text{GT}}\|_2$ \\ 
\midrule
$u_{\text{Tanh}}$ & 0.2 & 30 & 100 & 5,000 & 0.0071 \\ 
$u_{\text{ReLU}}$ & 0.2 & 30 & 100 & 5,000 & 0.0047 \\ 
$u_{\text{Step}}$  & 0.2 & 30 & 100 & 5,000 & 0.041 \\ 
$u_{\text{QStep}}$  & 0.25 & 150 & 100 & 8,000 & 0.067 \\ 
\bottomrule
\end{tabular}
\end{table}

\begin{figure}[htbp]
    \centering

    \begin{subfigure}{0.23\textwidth}
        \centering
        \includegraphics[width=\linewidth]{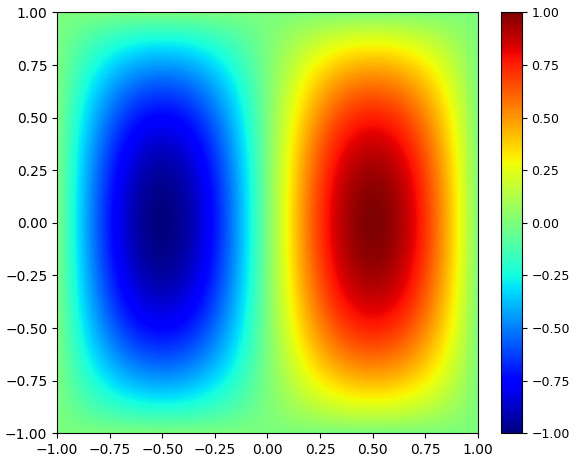}
        \caption{$u_{\text{GT}}$}
    \end{subfigure}\hfill
    \begin{subfigure}{0.23\textwidth}
        \centering
        \includegraphics[width=\linewidth]{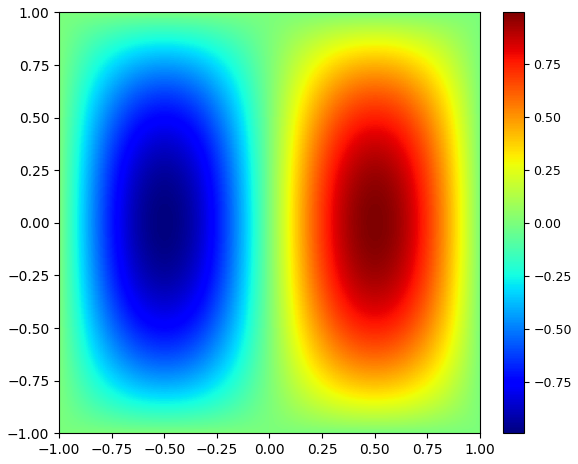}
        \caption{$u_{\text{Tanh}}$}
    \end{subfigure}\hfill
    \begin{subfigure}{0.24\textwidth}
        \centering
        \includegraphics[width=\linewidth]{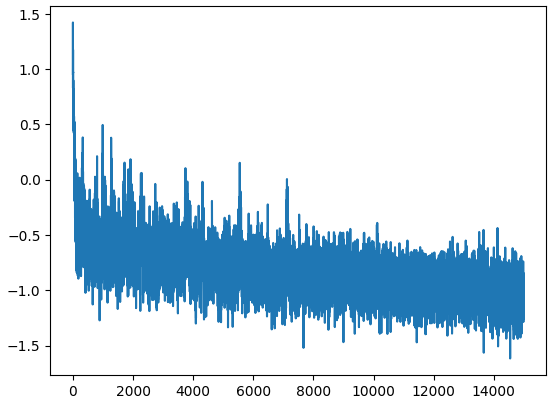}
        \caption{$\log(loss)$}
    \end{subfigure}\hfill
    \begin{subfigure}{0.24\textwidth}
        \centering
        \includegraphics[width=\linewidth]{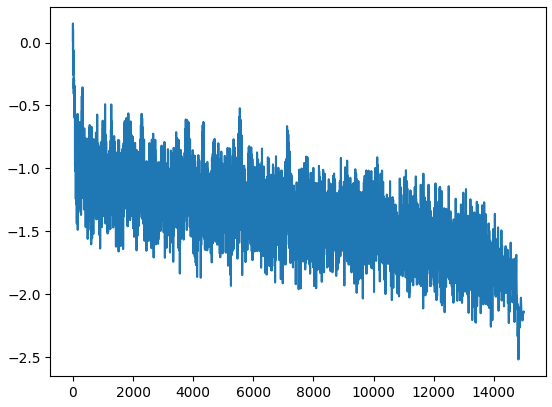}
        \caption{$\log( \|e_{\text{Tanh}}\|_{2} / \|u_{\text{GT}}\|_{2} )$}
    \end{subfigure}

    \vspace{1em} 

    \begin{subfigure}{0.23\textwidth}
        \centering
        \includegraphics[width=\linewidth]{suave_u_GT.png}
        \caption{$u_{\text{GT}}$}
    \end{subfigure}\hfill
    \begin{subfigure}{0.23\textwidth}
        \centering
        \includegraphics[width=\linewidth]{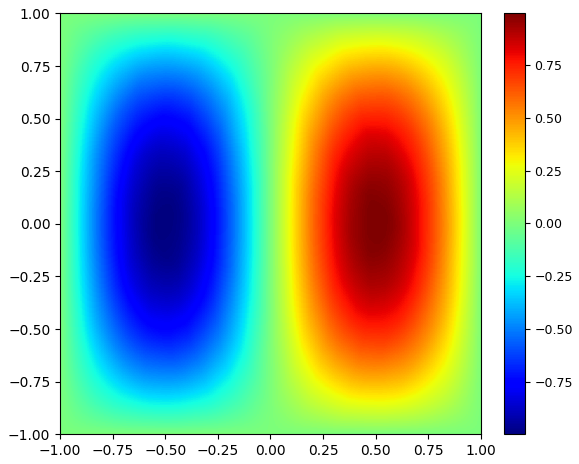}
        \caption{$u_{\text{ReLU}}$}
    \end{subfigure}\hfill
    \begin{subfigure}{0.24\textwidth}
        \centering
        \includegraphics[width=\linewidth]{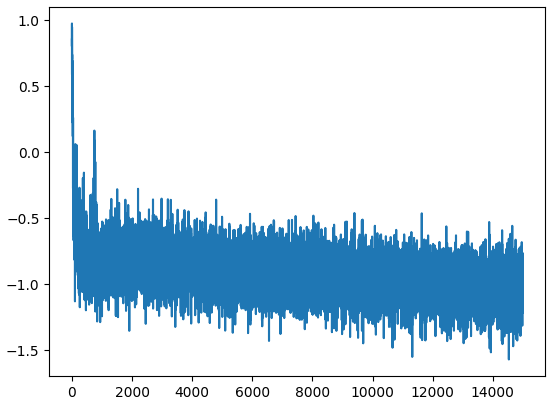}
        \caption{$\log(loss)$}
    \end{subfigure}\hfill
    \begin{subfigure}{0.24\textwidth}
        \centering
        \includegraphics[width=\linewidth]{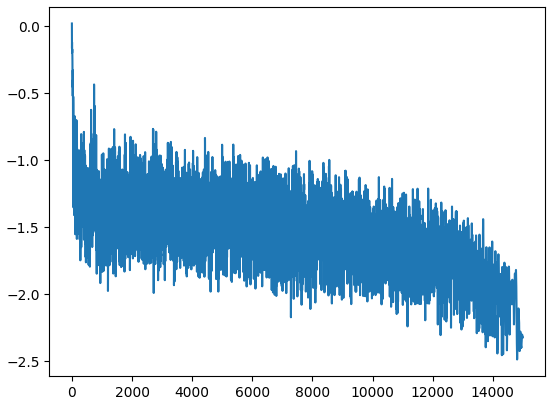}
        \caption{$\log(\|e_{\text{ReLU}}\|_{2} / \|u_{\text{GT}}\|_{2})$}
    \end{subfigure}

    \vspace{1em}

    \begin{subfigure}{0.23\textwidth}
        \centering
        \includegraphics[width=\linewidth]{suave_u_GT.png}
        \caption{$u_{\text{GT}}$}
    \end{subfigure}\hfill
    \begin{subfigure}{0.23\textwidth}
        \centering
        \includegraphics[width=\linewidth]{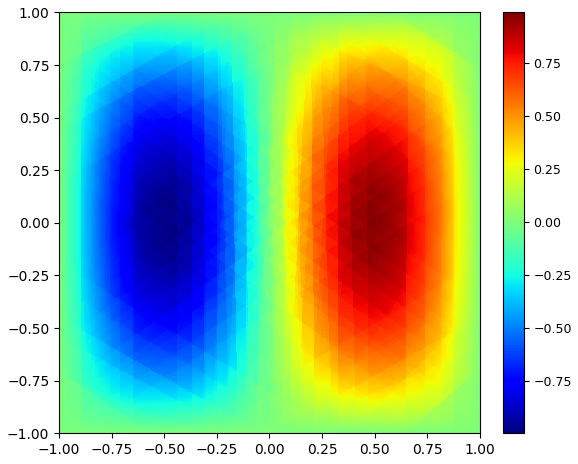}
        \caption{$u_{\text{Step}}$}
    \end{subfigure}\hfill
    \begin{subfigure}{0.24\textwidth}
        \centering
        \includegraphics[width=\linewidth]{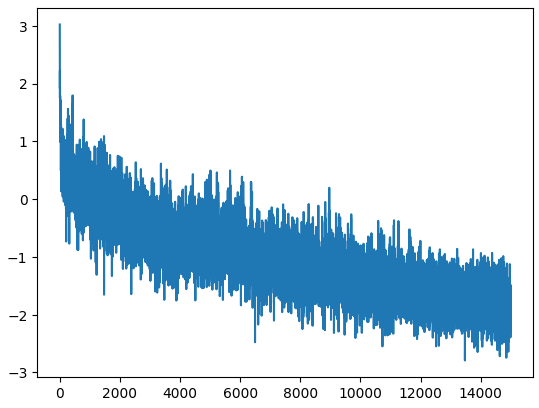}
        \caption{$\log(loss)$}
    \end{subfigure}\hfill
    \begin{subfigure}{0.24\textwidth}
        \centering
        \includegraphics[width=\linewidth]{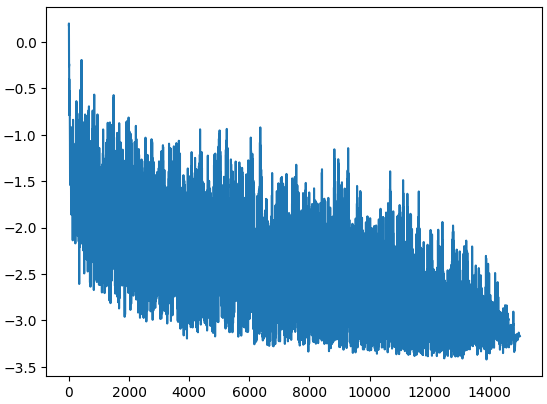}
        \caption{$\log(\|e_{\text{Step}}\|_{2} / \|u_{\text{GT}}\|_{2})$}
    \end{subfigure}

    \vspace{1em}

    \begin{subfigure}{0.23\textwidth}
        \centering
        \includegraphics[width=\linewidth]{suave_u_GT.png}
        \caption{$u_{\text{GT}}$}
    \end{subfigure}\hfill
    \begin{subfigure}{0.23\textwidth}
        \centering
        \includegraphics[width=\linewidth]{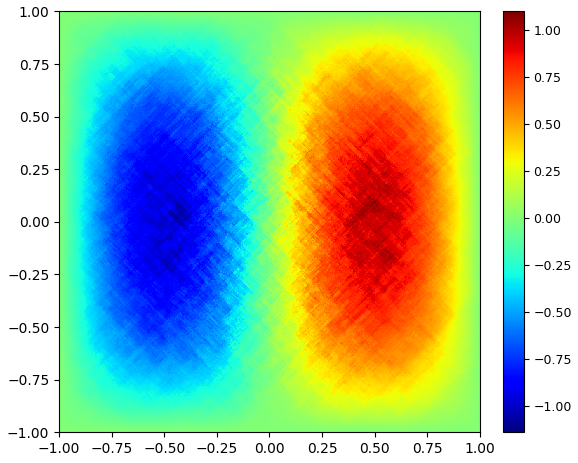}
        \caption{ $u_{\text{QStep}}$ }
    \end{subfigure}\hfill
    \begin{subfigure}{0.24\textwidth}
        \centering
        \includegraphics[width=\linewidth]{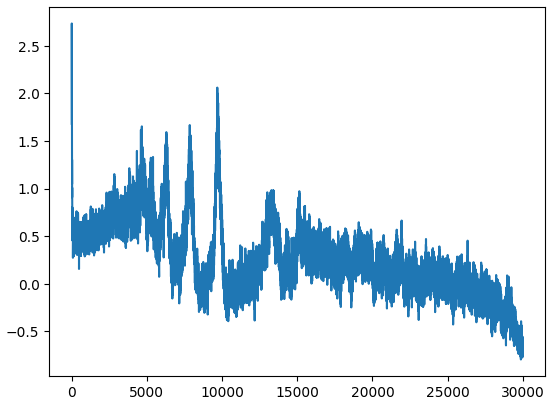}
        \caption{$\log{(loss)}$}
    \end{subfigure}\hfill
    \begin{subfigure}{0.24\textwidth}
        \centering
        \includegraphics[width=\linewidth]{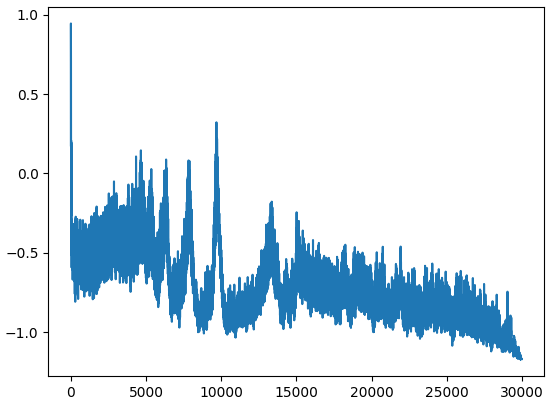}
        \caption{$\log( \|e_{\text{QStep}}\|_{2} / \|u_{\text{GT}}\|_{2} )$}
    \end{subfigure}

    \caption{Numerical results for Example \ref{sec:suave} using Algorithm \ref{alg:flujo}. Rows correspond to architectures. Columns (left to right): exact solution, predicted solution, loss history (log scale), and relative $L^2$ error history (log scale).}
    \label{fig:suave}
\end{figure}

\end{example}

\begin{example}[High-dimensional case]\label{sec:hdim}

We consider nonlinear second-order Sine-Gordon equations characterized by anisotropic and inseparable solutions. To evaluate the performance of the proposed method in high-dimensional regimes, we adopt the following exact solution established in the literature \cite{hdim1, hdim2, hdim3}:

\begin{equation}
u_{\text{GT}}(x) = \left(1 - \|x\|^2\right) \left( \sum_{i=1}^{d-1} c_i \sin(x_i + \cos(x_{i+1}) + x_{i+1} \cos(x_i)) \right),
\end{equation}
where the coefficients $c_i$ are sampled from a standard normal distribution $\mathcal{N}(0, 1)$. This solution exhibits correlation between the $i$-th and $(i+1)$-th coordinates, frequently referred to as two-body interaction.

Given this exact solution, we consider the Sine-Gordon equation subject to homogeneous Dirichlet boundary conditions:

\begin{equation}\label{eq:hdim}
\left\lbrace \begin{aligned}
	\Delta u(x) + \sin (u(x)) &= g(x), \quad & \mbox{in } \Omega,\\
	u & = 0 \quad & \mbox{on } \partial\Omega,\\
\end{aligned} \right.
\end{equation}
where $\Omega = \{x \in \R^d : \|x\| < 1\}$ and the source term is prescribed as $g(x) = \Delta u_{\text{GT}}(x) + \sin (u_{\text{GT}}(x))$.

To approximate the solution in this test case, we implement the framework described in Algorithm~\ref{alg:dominio_ficticio}. Following the architectural choices of our previous experiment, all trial functions are parameterized using multilayer neural networks equipped with two hidden layers of $128$ neurons each. In particular, the formulation for $u_{\text{Step}}$ involves a network mapping characterized by the configuration $\mathcal{N}(4, (128,128), 4)$, as detailed in Section~\ref{subsec:step_arch}. The underlying system relies on a collection of randomly initialized coefficients $c_i$, which are explicitly given by:
\begin{align*}
    c = [&-0.2793, -1.5054, \phantom{-}0.6111, \phantom{-}0.7158, -0.9995, \phantom{-}1.2583, -1.1201, \phantom{-}1.1827, \\
    &\phantom{-}0.1987, -1.1433, -1.3718, -0.3441, -1.0641, \phantom{-}0.7060, -2.0442, \phantom{-}0.1833, \\
    &\phantom{-}0.0275, -0.1868, \phantom{-}1.0890]
\end{align*}
The resulting numerical performance, alongside the comprehensive hyperparameter configurations, is documented in Table~\ref{table:hdim} and visualized in Figure~\ref{fig:hdim}.

\begin{table}[htbp]
\centering
\caption{Hyperparameters and resulting relative errors for Algorithm \ref{alg:dominio_ficticio} applied to Example \ref{sec:hdim} (High-dimensional case).}
\label{table:hdim}
\begin{tabular}{l ccc r}
\toprule
Architecture & $h$ & $N_x$ & $N_y$ & $\|e_*\|_2 / \|u_{\text{GT}}\|_2$ \\ 
\midrule
$u_{\text{Tanh}}$ & 0.15 & 400 & 100 & 0.011 \\ 
$u_{\text{ReLU}}$ & 0.15 & 400 & 100 & 0.013 \\ 
$u_{\text{Step}}$  & 0.15 & 400 & 100 & 0.020 \\ 
\bottomrule
\end{tabular}
\end{table}

\begin{figure}[htbp]
    \centering

    \begin{subfigure}{0.22\textwidth}
        \centering
        \includegraphics[width=\linewidth]{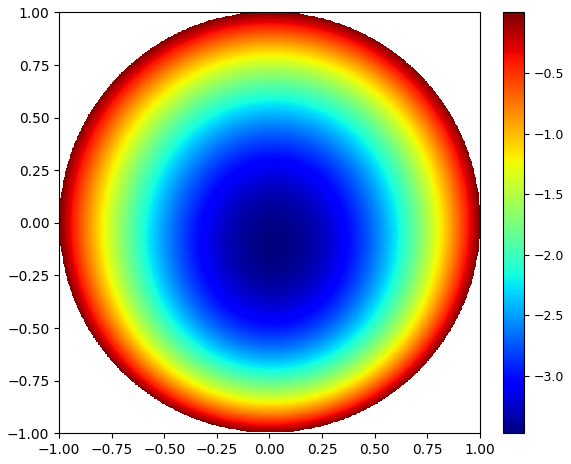}
        \caption{$u_{\text{GT}}$}
    \end{subfigure}\hfill
    \begin{subfigure}{0.22\textwidth}
        \centering
        \includegraphics[width=\linewidth]{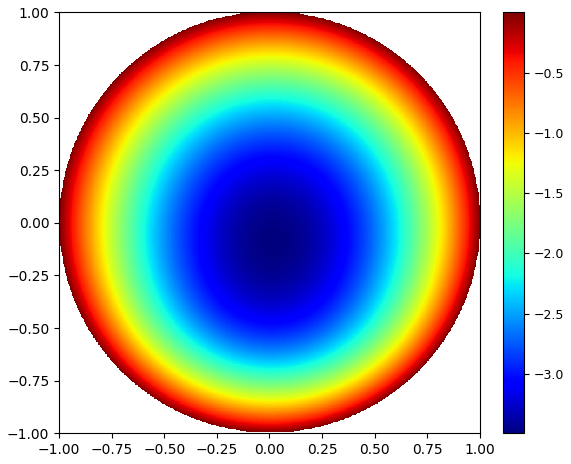}
        \caption{$u_{\text{Tanh}}$}
    \end{subfigure}\hfill
    \begin{subfigure}{0.24\textwidth}
        \centering
        \includegraphics[width=\linewidth]{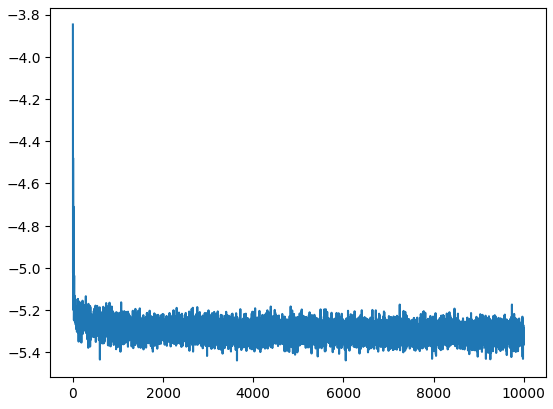}
        \caption{$\log(loss)$}
    \end{subfigure}\hfill
    \begin{subfigure}{0.24\textwidth}
        \centering
        \includegraphics[width=\linewidth]{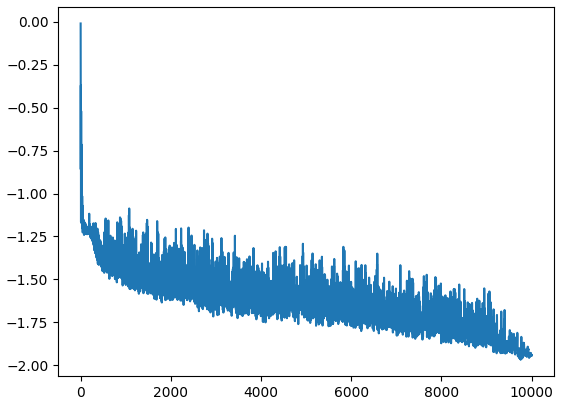}
        \caption{$\log( \|e_{\text{Tanh}}\|_{2} / \|u_{\text{GT}}\|_{2} )$}
    \end{subfigure}

    \vspace{1em} 

    \begin{subfigure}{0.22\textwidth}
        \centering
        \includegraphics[width=\linewidth]{20dim_u_GT.png}
        \caption{$u_{\text{GT}}$}
    \end{subfigure}\hfill
    \begin{subfigure}{0.22\textwidth}
        \centering
        \includegraphics[width=\linewidth]{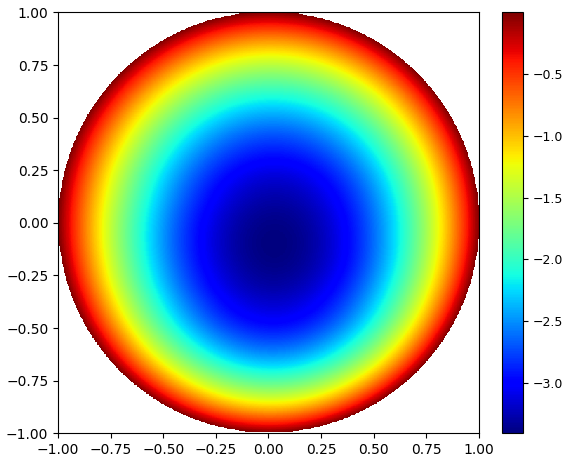}
        \caption{$u_{\text{ReLU}}$}
    \end{subfigure}\hfill
    \begin{subfigure}{0.24\textwidth}
        \centering
        \includegraphics[width=\linewidth]{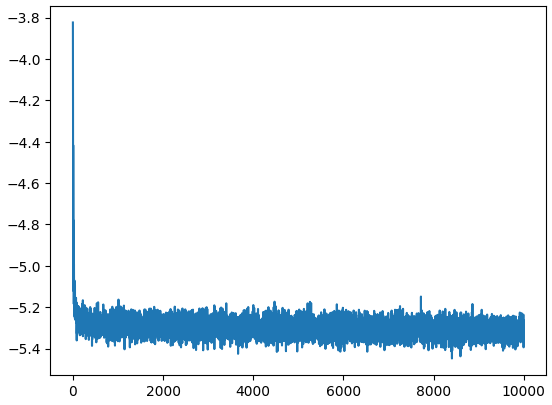}
        \caption{$\log(loss)$}
    \end{subfigure}\hfill
    \begin{subfigure}{0.24\textwidth}
        \centering
        \includegraphics[width=\linewidth]{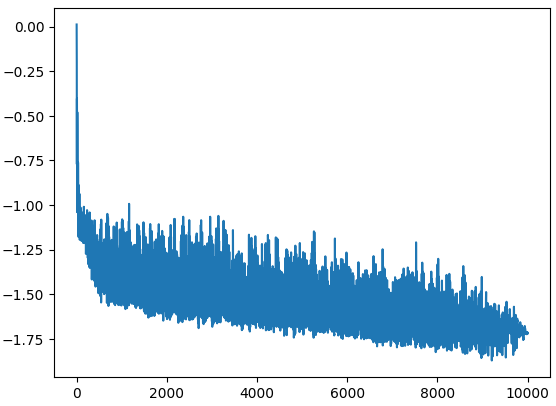}
        \caption{$\log(\|e_{\text{ReLU}}\|_{2} / \|u_{\text{GT}}\|_{2})$}
    \end{subfigure}

    \vspace{1em}

    \begin{subfigure}{0.22\textwidth}
        \centering
        \includegraphics[width=\linewidth]{20dim_u_GT.png}
        \caption{$u_{\text{GT}}$}
    \end{subfigure}\hfill
    \begin{subfigure}{0.22\textwidth}
        \centering
        \includegraphics[width=\linewidth]{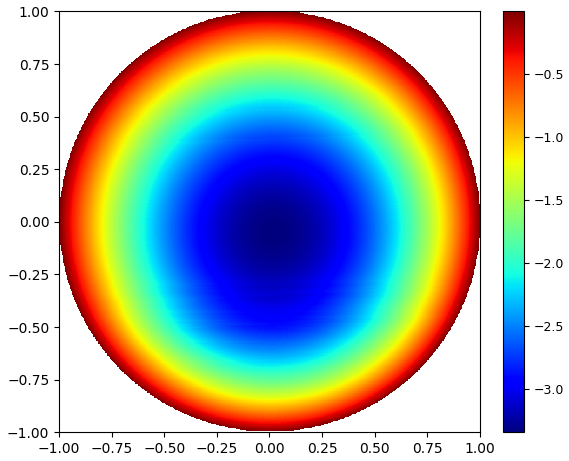}
        \caption{$u_{\text{Step}}$}
    \end{subfigure}\hfill
    \begin{subfigure}{0.24\textwidth}
        \centering
        \includegraphics[width=\linewidth]{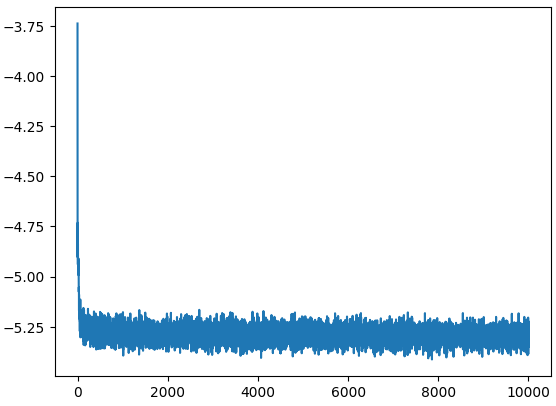}
        \caption{$\log(loss)$}
    \end{subfigure}\hfill
    \begin{subfigure}{0.24\textwidth}
        \centering
        \includegraphics[width=\linewidth]{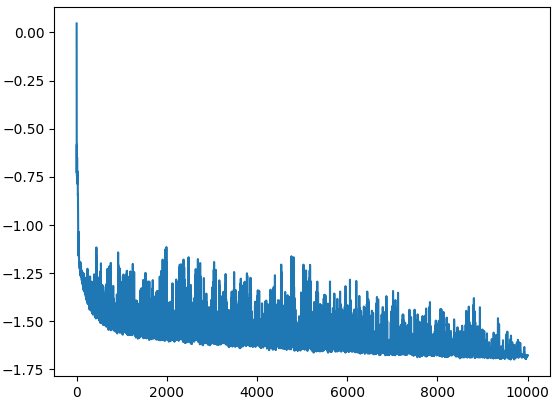}
        \caption{$\log(\|e_{\text{Step}}\|_{2} / \|u_{\text{GT}}\|_{2})$}
    \end{subfigure}

    \vspace{1em}

    \caption{
Numerical results for Example \ref{sec:hdim} using Algorithm \ref{alg:dominio_ficticio}. Rows correspond to architectures. Columns (left to right): 2D slice of the exact solution, 2D slice of the predicted solution, loss history (log scale), and relative $L^2$ error history (log scale).
    }
    \label{fig:hdim}
\end{figure}

\end{example}

\begin{example}[Point source term]\label{sec:dirac}

To further evaluate the proposed method, we consider the Poisson equation with a Dirac delta distribution supported at a single point. Let $\Omega = B(0,1) \subset \mathbb{R}^2$ denote the unit disk centered at the origin. We seek a solution $u$ satisfying
\begin{equation}
\begin{cases} 
-\Delta u = \delta & \text{in } \Omega, \\ 
u = 0 & \text{on } \partial\Omega,
\end{cases}
\end{equation}
where $\delta$ denotes the Dirac measure supported at the origin. The exact solution is given by the fundamental solution of the Laplacian in two dimensions, modified to satisfy the homogeneous Dirichlet boundary condition:
\begin{equation}
u_{\text{GT}}(x) = -\frac{1}{2\pi} \log(\|x\|).
\end{equation}

We approximate the solution using Algorithm \ref{alg:dominio_ficticio}. As before, the trial function approximation spaces employ neural networks with two $128$-neuron hidden layers. Specifically,  $u_{\text{Step}}$ uses the topology $\mathcal{N}(6, (64, 64, 64, 64), 16)$ (Section \ref{subsec:step_arch}). Table \ref{table:dirac} and Figure \ref{fig:dirac} report the hyperparameters and numerical results.

In this example, the non-smooth architectures $u_{\text{ReLU}}$ and $u_{\text{Step}}$ outperform $u_{\text{Tanh}}$. This performance gap is directly attributable to the singularity in the solution, a regime where non-smooth architectures inherently offer superior approximation capabilities.

\begin{table}[htbp]
\centering
\caption{Hyperparameters and resulting relative errors for Algorithm \ref{alg:dominio_ficticio} applied to Example \ref{sec:dirac} (Point source term).}
\label{table:dirac}
\begin{tabular}{l ccc r}
\toprule
Architecture & $h$ & $N_x$ & $N_y$ & $\|e_*\|_2 / \|u_{\text{GT}}\|_2$ \\ 
\midrule
$u_{\text{Tanh}}$ & 0.2 & 25 & 200 & 0.023 \\ 
$u_{\text{ReLU}}$ & 0.2 & 25 & 200 & 0.013 \\ 
$u_{\text{Step}}$  & 0.2 & 200 & 20 & 0.047 \\ 
\bottomrule
\end{tabular}
\end{table}

\begin{figure}[htbp]
    \centering

    \begin{subfigure}{0.22\textwidth}
        \centering
        \includegraphics[width=\linewidth]{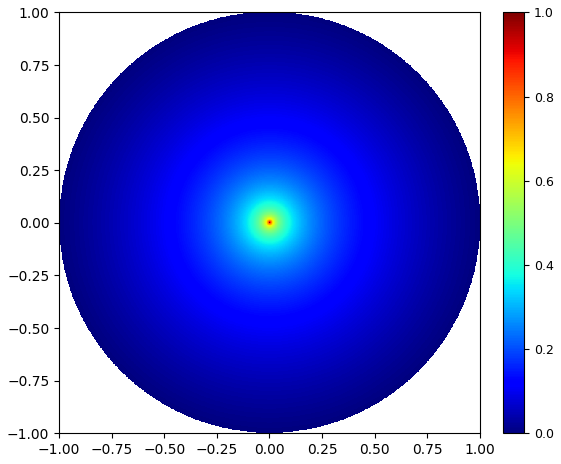}
        \caption{$u_{\text{GT}}$}
    \end{subfigure}\hfill
    \begin{subfigure}{0.22\textwidth}
        \centering
        \includegraphics[width=\linewidth]{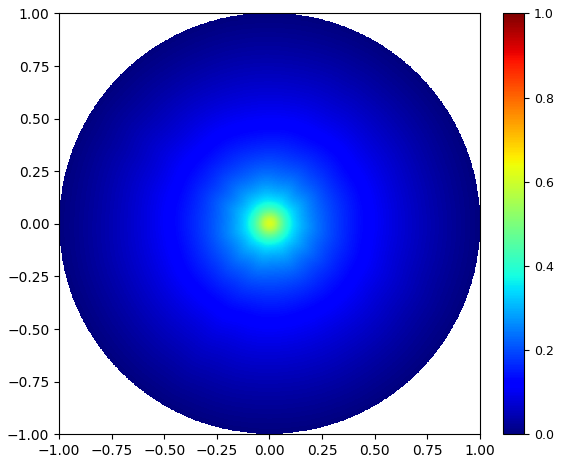}
        \caption{$u_{\text{Tanh}}$}
    \end{subfigure}\hfill
    \begin{subfigure}{0.24\textwidth}
        \centering
        \includegraphics[width=\linewidth]{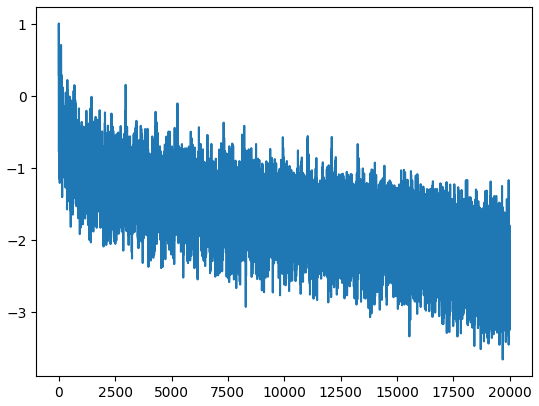}
        \caption{$\log(loss)$}
    \end{subfigure}\hfill
    \begin{subfigure}{0.24\textwidth}
        \centering
        \includegraphics[width=\linewidth]{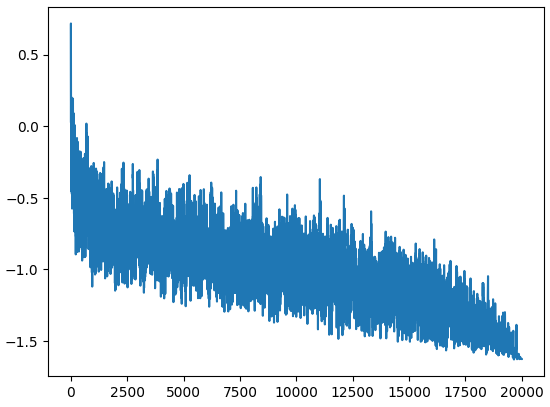}
        \caption{$\log( \|e_{\text{Tanh}}\|_{2} / \|u_{\text{GT}}\|_{2} )$}
    \end{subfigure}

    \vspace{1em} 

    \begin{subfigure}{0.22\textwidth}
        \centering
        \includegraphics[width=\linewidth]{dirac_u_GT.png}
        \caption{$u_{\text{GT}}$}
    \end{subfigure}\hfill
    \begin{subfigure}{0.22\textwidth}
        \centering
        \includegraphics[width=\linewidth]{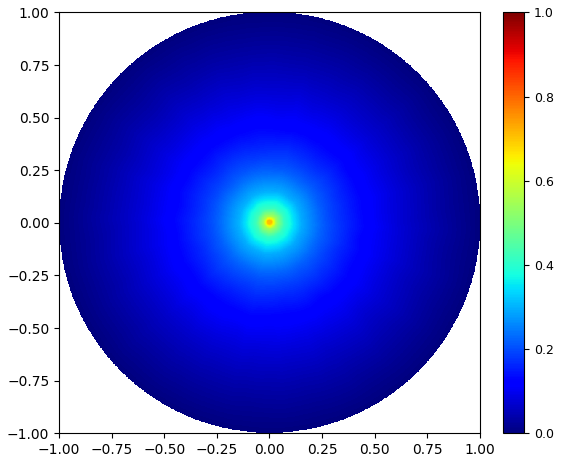}
        \caption{$u_{\text{ReLU}}$}
    \end{subfigure}\hfill
    \begin{subfigure}{0.24\textwidth}
        \centering
        \includegraphics[width=\linewidth]{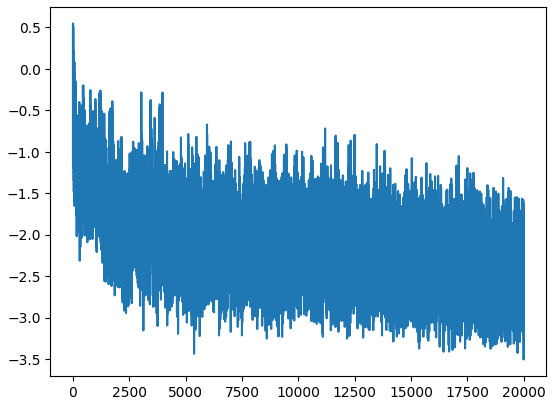}
        \caption{$\log(loss)$}
    \end{subfigure}\hfill
    \begin{subfigure}{0.24\textwidth}
        \centering
        \includegraphics[width=\linewidth]{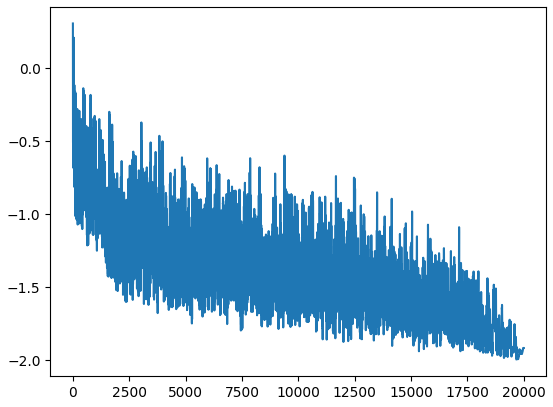}
        \caption{$\log(\|e_{\text{ReLU}}\|_{2} / \|u_{\text{GT}}\|_{2})$}
    \end{subfigure}

    \vspace{1em}

    \begin{subfigure}{0.22\textwidth}
        \centering
        \includegraphics[width=\linewidth]{dirac_u_GT.png}
        \caption{$u_{\text{GT}}$}
    \end{subfigure}\hfill
    \begin{subfigure}{0.22\textwidth}
        \centering
        \includegraphics[width=\linewidth]{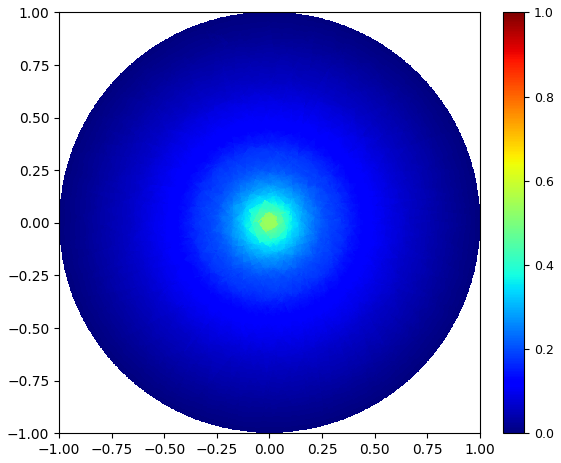}
        \caption{$u_{\text{Step}}$}
    \end{subfigure}\hfill
    \begin{subfigure}{0.24\textwidth}
        \centering
        \includegraphics[width=\linewidth]{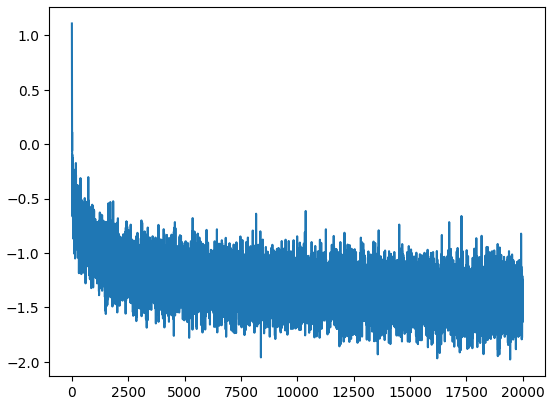}
        \caption{$\log(loss)$}
    \end{subfigure}\hfill
    \begin{subfigure}{0.24\textwidth}
        \centering
        \includegraphics[width=\linewidth]{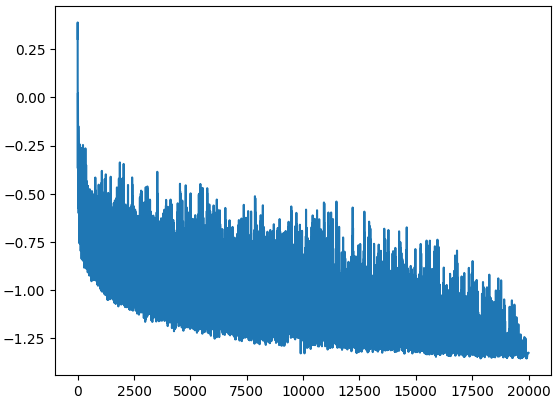}
        \caption{$\log(\|e_{\text{Step}}\|_{2} / \|u_{\text{GT}}\|_{2})$}
    \end{subfigure}

    \vspace{1em}

    \caption{
    Results for Example \ref{sec:dirac} using Algorithm \ref{alg:dominio_ficticio}. Rows correspond to the three architectural setups. Columns display the exact solution, approximate solution, log-scale loss history, and logarithmic relative $L^2$ error versus training iterations.
    }
    \label{fig:dirac}
\end{figure}

\end{example}

\begin{example}[Segment source term]\label{sec:segment}

We consider a Poisson equation featuring a singular source term supported on a one-dimensional segment. Let $\xi$ denote the segment defined by
\[
\xi = \{(0, 0, z) \in \mathbb{R}^3 : -1 \le z \le 1\}.
\]
We prescribe the source term $\gamma = \delta_\xi$. The corresponding free-space fundamental solution is expressed by the integral
\begin{align*}
G(x,y,z) &= \int_{-1}^{1} \frac{1}{4\pi\sqrt{x^2+y^2+(z-t)^2}} \, dt \\
&= \frac{1}{4\pi} \log \left( \frac{\sqrt{x^2+y^2+(z-1)^2} + 1 - z}{\sqrt{x^2+y^2+(z+1)^2} - 1 - z} \right).
\end{align*}

For our numerical experiments, we define the domain $\Omega$ as the solid ellipsoid bounded by the level set where the argument of the logarithm in $G$ equals $3$:
\[
\Omega = \left\{ (x, y, z) \in \mathbb{R}^3 : \frac{x^2}{3} + \frac{y^2}{3} + \frac{z^2}{4} \le 1 \right\}.
\]

The boundary value problem seeks a function $u$ such that
\begin{equation}
\begin{cases} 
-\Delta u = \delta_\xi & \text{in } \Omega, \\ 
u = 0 & \text{on } \partial\Omega. 
\end{cases}
\end{equation}

The exact singular solution to this problem is given by
\begin{equation}
u_{\text{GT}}(x,y,z) = \frac{1}{4\pi} \log \left( \frac{\sqrt{x^2+y^2+(z-1)^2} + 1 - z}{3(\sqrt{x^2+y^2+(z+1)^2} - 1 - z)} \right).
\end{equation}

In this example, we approximate the solution using Algorithm \ref{alg:dominio_ficticio} and two-layer neural networks with $128$ neurons each. Specifically, the architecture for $u_{\text{Step}}$ is $\mathcal{N}(6, (64, 64, 64, 64), 16)$ (Section \ref{subsec:step_arch}). Table \ref{table:segment} and Figure \ref{fig:segment} report the hyperparameters and numerical results.

As in the previous example, the non-smooth architectures $u_{\text{ReLU}}$ and $u_{\text{Step}}$ exhibit superior overall performance compared to $u_{\text{Tanh}}$. Consistently, non-smooth architectures prove more suitable for approximating singular solutions.

\begin{table}[htbp]
\centering
\caption{Hyperparameters and resulting relative errors for Algorithm \ref{alg:dominio_ficticio} applied to Example \ref{sec:segment} (Segment source term).}
\label{table:segment}
\begin{tabular}{l ccc r}
\toprule
Architecture & $h$ & $N_x$ & $N_y$ & $\|e_*\|_2 / \|u_{\text{GT}}\|_2$ \\ 
\midrule
$u_{\text{Tanh}}$ & 0.4 & 200 & 40 & 0.059 \\ 
$u_{\text{ReLU}}$ & 0.4 & 200 & 40 & 0.047 \\ 
$u_{\text{Step}}$  & 0.4 & 400 & 40 & 0.066 \\ 
\bottomrule
\end{tabular}
\end{table}

\begin{figure}[htbp]
    \centering
    
    \begin{subfigure}{0.24\textwidth}
        \centering
        \includegraphics[width=\linewidth]{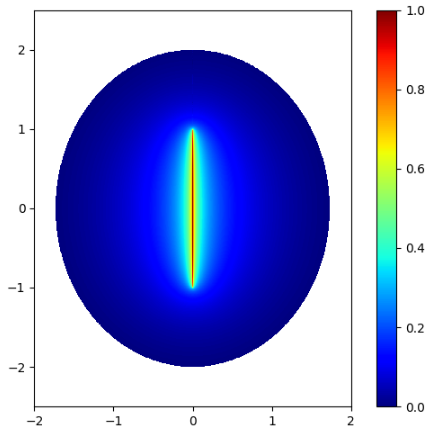}
        \caption{$u_{\text{GT}}$}
    \end{subfigure}\hfill
    \begin{subfigure}{0.24\textwidth}
        \centering
        \includegraphics[width=\linewidth]{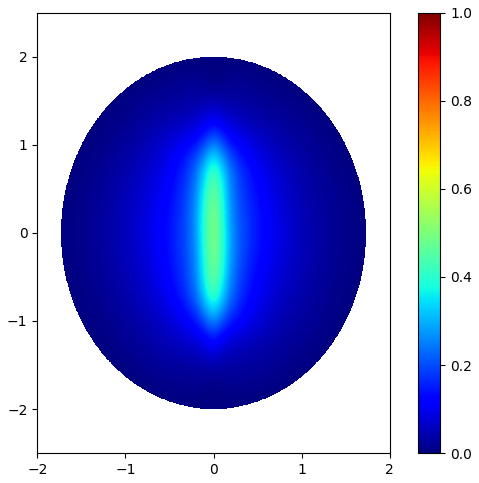}
        \caption{$u_{\text{Tanh}}$}
    \end{subfigure}\hfill
    \begin{subfigure}{0.24\textwidth}
        \centering
        \includegraphics[width=\linewidth]{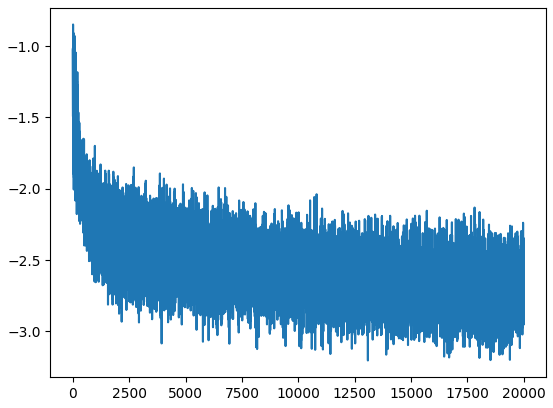}
        \caption{$\log(loss)$}
    \end{subfigure}\hfill
    \begin{subfigure}{0.24\textwidth}
        \centering
        \includegraphics[width=\linewidth]{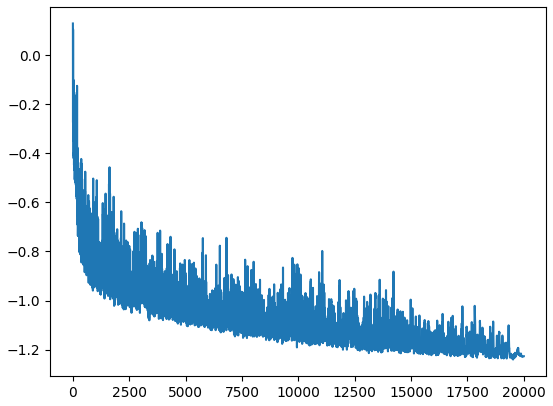}
        \caption{$\log( \|e_{\text{Tanh}}\|_{2} / \|u_{\text{GT}}\|_{2} )$}
    \end{subfigure}

    \vspace{1em} 

    \begin{subfigure}{0.24\textwidth}
        \centering
        \includegraphics[width=\linewidth]{segment_u_GT.png}
        \caption{$u_{\text{GT}}$}
    \end{subfigure}\hfill
    \begin{subfigure}{0.24\textwidth}
        \centering
        \includegraphics[width=\linewidth]{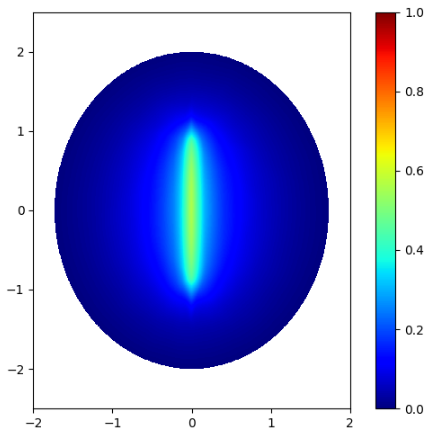}
        \caption{$u_{\text{ReLU}}$}
    \end{subfigure}\hfill
    \begin{subfigure}{0.24\textwidth}
        \centering
        \includegraphics[width=\linewidth]{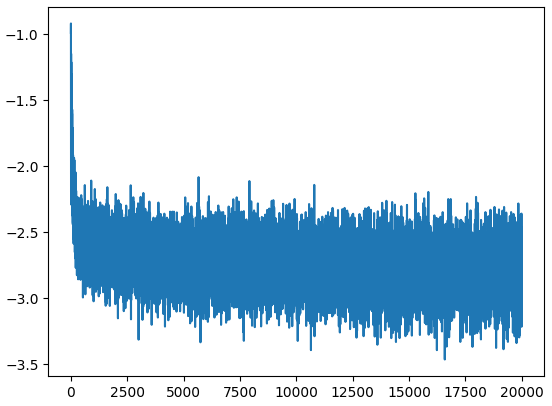}
        \caption{$\log(loss)$}
    \end{subfigure}\hfill
    \begin{subfigure}{0.24\textwidth}
        \centering
        \includegraphics[width=\linewidth]{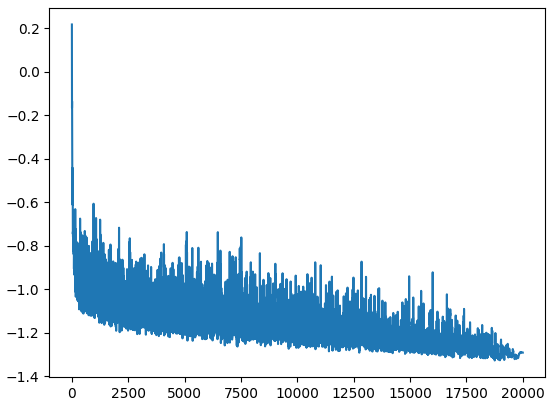}
        \caption{$\log(\|e_{\text{ReLU}}\|_{2} / \|u_{\text{GT}}\|_{2})$}
    \end{subfigure}

    \vspace{1em}

    \begin{subfigure}{0.24\textwidth}
        \centering
        \includegraphics[width=\linewidth]{segment_u_GT.png}
        \caption{$u_{\text{GT}}$}
    \end{subfigure}\hfill
    \begin{subfigure}{0.24\textwidth}
        \centering
        \includegraphics[width=\linewidth]{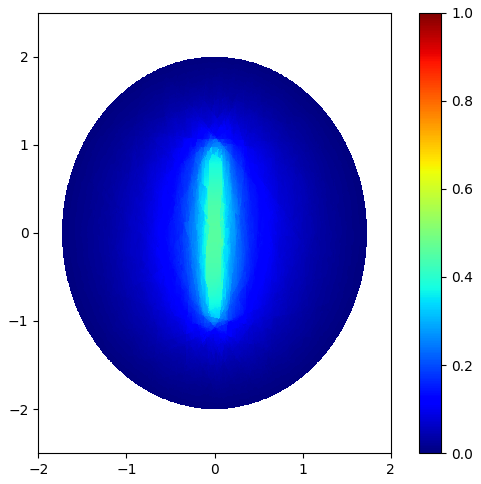}
        \caption{$u_{\text{Step}}$}
    \end{subfigure}\hfill
    \begin{subfigure}{0.24\textwidth}
        \centering
        \includegraphics[width=\linewidth]{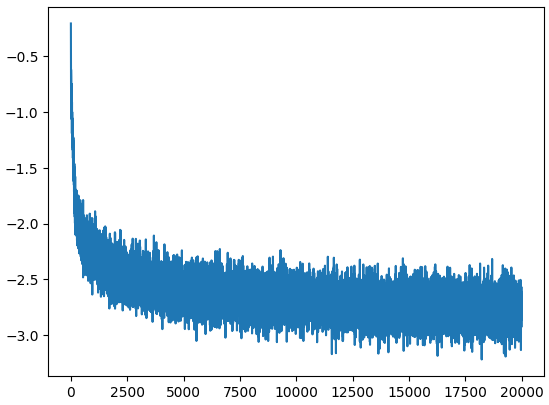}
        \caption{$\log(loss)$}
    \end{subfigure}\hfill
    \begin{subfigure}{0.24\textwidth}
        \centering
        \includegraphics[width=\linewidth]{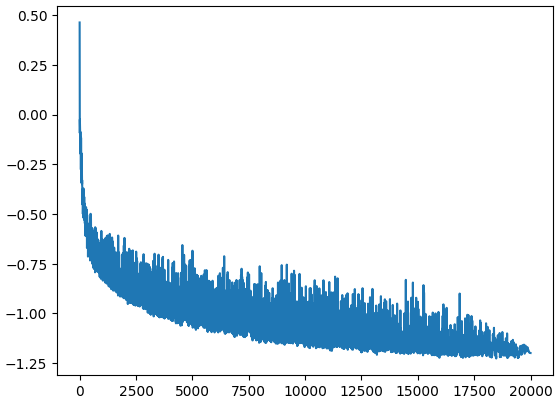}
        \caption{$\log(\|e_{\text{Step}}\|_{2} / \|u_{\text{GT}}\|_{2})$}
    \end{subfigure}

    \vspace{1em}

    \caption{
    Results for Example \ref{sec:segment} using Algorithm \ref{alg:dominio_ficticio}. Rows correspond to three neural architectures. Columns display 2D slices of the exact and approximate solutions, log-scale loss history, and logarithmic relative $L^2$ error.
    }
    \label{fig:segment}
\end{figure}

\end{example}


\begin{example}[L-shaped domain]\label{sec:pacman}

We now consider a non-smooth domain featuring a re-entrant corner, a configuration where the solution is expected to exhibit a singular derivative. Let $\Omega = B_1(0) \setminus [0,1]^2$. We seek a function $u: \Omega \to \mathbb{R}$ such that
\begin{equation}
\begin{cases} 
-\Delta u = 1 & \text{in } \Omega, \\ 
u = 0 & \text{on } \partial\Omega,
\end{cases}
\end{equation}
In the absence of a closed-form expression for the exact solution, we evaluate the performance of our method against a reference solution, denoted by $u_{\text{FEM}}$, obtained via the finite element method. The reference solution is computed using linear finite elements on a triangular mesh of approximately 20,000 nodes, with a mesh grading toward the origin to resolve the local singularity.

To accommodate the non-smooth boundary $\partial\Omega$, we apply Algorithm~\ref{alg:flujo}. All trial functions are modeled using neural networks with two $128$-neuron hidden layers. Specifically, the architecture for $u_{\text{Step}}$ is $\mathcal{N}(4, (64, 64, 64, 64), 16)$ (Section \ref{subsec:step_arch}). Table \ref{table:pacman} and Figure \ref{fig:pacman} report the hyperparameters and numerical results.

Unlike the two preceding examples, the solution in this case exhibits a singularity only in its derivatives. Nonetheless, the smooth architecture $u_{\text{Tanh}}$ underperforms relative to its non-smooth counterparts, $u_{\text{ReLU}}$ and $u_{\text{Step}}$.

\begin{table}[htbp]
\centering
\caption{Hyperparameters and resulting relative errors for Algorithm \ref{alg:flujo} applied to Example \ref{sec:pacman} (L-shaped domain).}
\label{table:pacman}
\begin{tabular}{l cccc r}
\toprule
Architecture & $h$ & $N_x$ & $N_y$ & $N_b$ & $\|e_*\|_2 / \|u_{\text{GT}}\|_2$ \\ 
\midrule
$u_{\text{Tanh}}$ & 0.2 & 40 & 100 & 3,000 & 0.031 \\ 
$u_{\text{ReLU}}$ & 0.2 & 40 & 100 & 3,000 & 0.017 \\ 
$u_{\text{Step}}$  & 0.2 & 100 & 40 & 3,000 & 0.030 \\ 
\bottomrule
\end{tabular}
\end{table}

\begin{figure}[htbp]
    \centering

    \begin{subfigure}{0.24\textwidth}
        \centering
        \includegraphics[width=\linewidth]{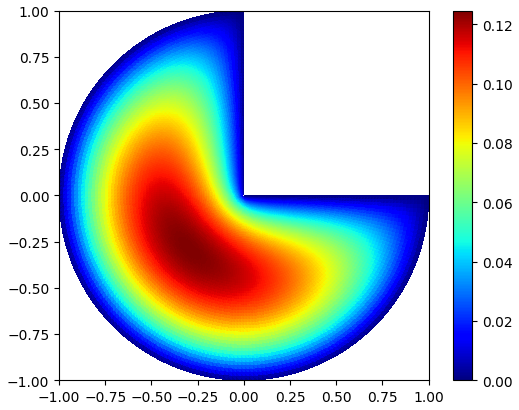}
        \caption{$u_{\text{FEM}}$}
    \end{subfigure}\hfill
    \begin{subfigure}{0.24\textwidth}
        \centering
        \includegraphics[width=\linewidth, height=2.85cm]{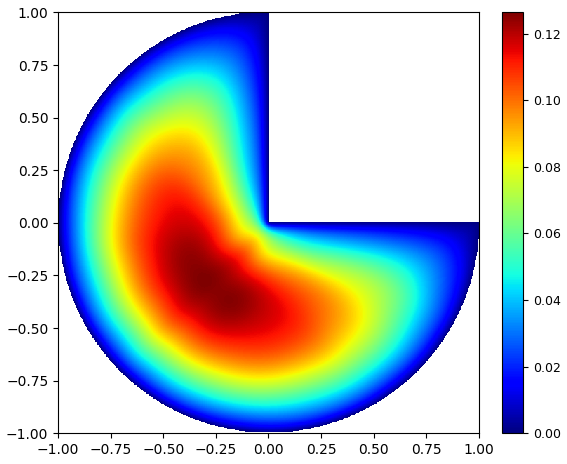}
        \caption{$u_{\text{Tanh}}$}
    \end{subfigure}\hfill
    \begin{subfigure}{0.24\textwidth}
        \centering
        \includegraphics[width=\linewidth]{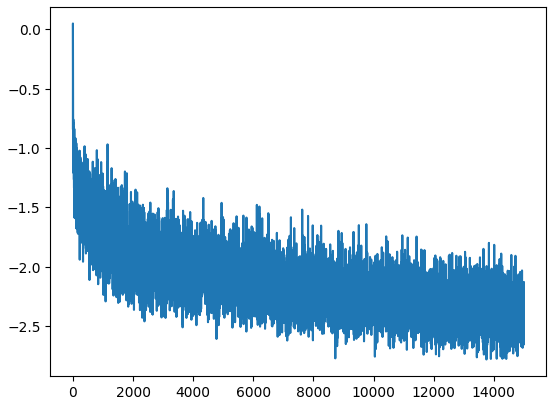}
        \caption{$\log(loss)$}
    \end{subfigure}\hfill
    \begin{subfigure}{0.24\textwidth}
        \centering
        \includegraphics[width=\linewidth]{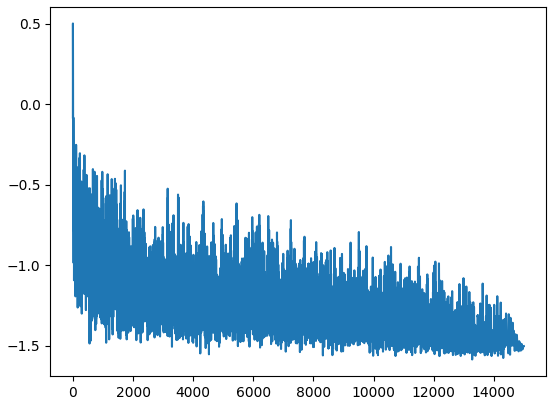}
        \caption{$\log( \|e_{\text{Tanh}}\|_{2} / \|u_{\text{GT}}\|_{2} )$}
    \end{subfigure}

    \vspace{1em} 

    \begin{subfigure}{0.24\textwidth}
        \centering
        \includegraphics[width=\linewidth]{pacman_u_GT.png}
        \caption{$u_{\text{FEM}}$}
    \end{subfigure}\hfill
    \begin{subfigure}{0.24\textwidth}
        \centering
        \includegraphics[width=\linewidth, height=2.85cm]{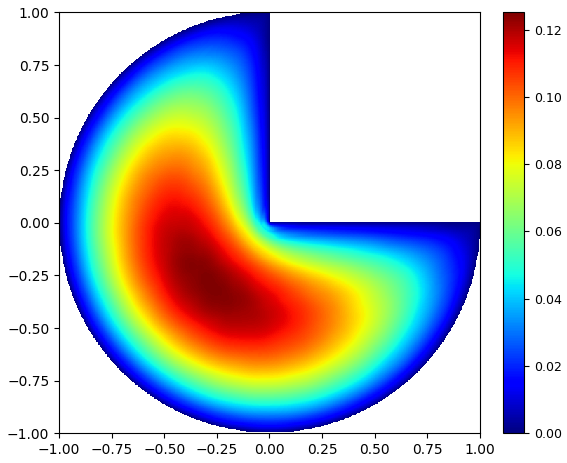}
        \caption{$u_{\text{ReLU}}$}
    \end{subfigure}\hfill
    \begin{subfigure}{0.24\textwidth}
        \centering
        \includegraphics[width=\linewidth]{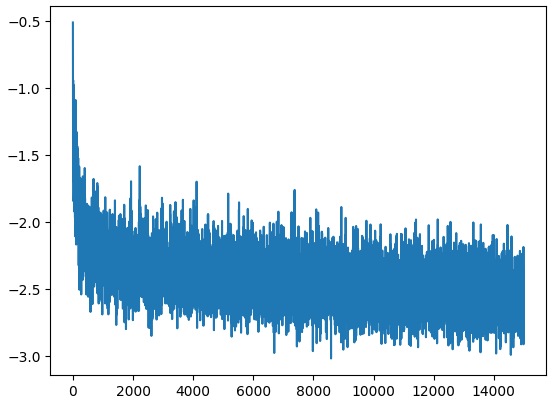}
        \caption{$\log(loss)$}
    \end{subfigure}\hfill
    \begin{subfigure}{0.24\textwidth}
        \centering
        \includegraphics[width=\linewidth]{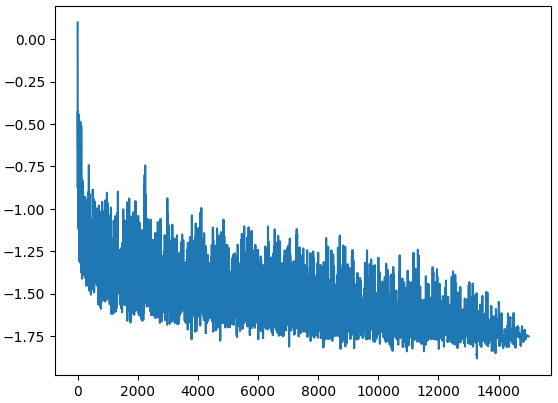}
        \caption{$\log(\|e_{\text{ReLU}}\|_{2} / \|u_{\text{GT}}\|_{2})$}
    \end{subfigure}

    \vspace{1em}

    \begin{subfigure}{0.24\textwidth}
        \centering
        \includegraphics[width=\linewidth]{pacman_u_GT.png}
        \caption{$u_{\text{FEM}}$}
    \end{subfigure}\hfill
    \begin{subfigure}{0.24\textwidth}
        \centering
        \includegraphics[width=\linewidth, height=2.85cm]{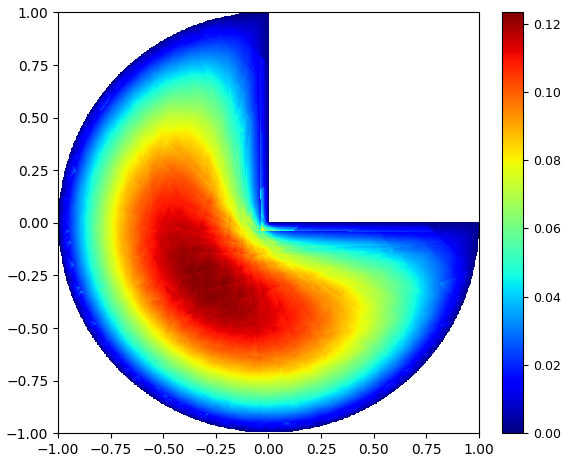}
        \caption{$u_{\text{Step}}$}
    \end{subfigure}\hfill
    \begin{subfigure}{0.24\textwidth}
        \centering
        \includegraphics[width=\linewidth]{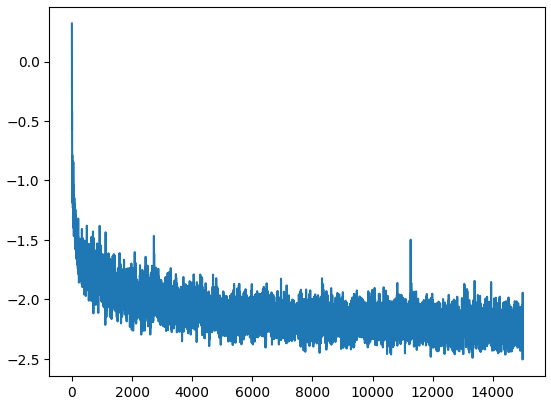}
        \caption{$\log(loss)$}
    \end{subfigure}\hfill
    \begin{subfigure}{0.24\textwidth}
        \centering
        \includegraphics[width=\linewidth]{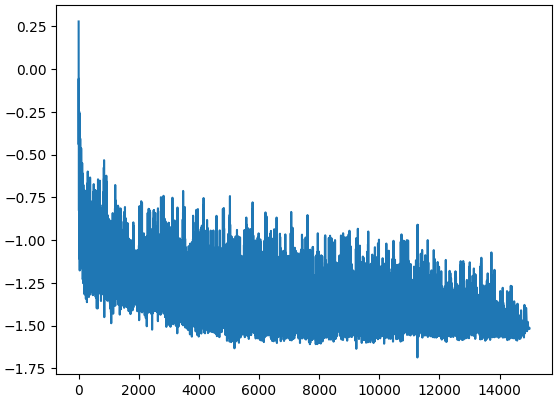}
        \caption{$\log(\|e_{\text{Step}}\|_{2} / \|u_{\text{GT}}\|_{2})$}
    \end{subfigure}

    \vspace{1em}

    \caption{
   Results for Example \ref{sec:pacman} using Algorithm \ref{alg:flujo}. Rows correspond to three network architectures. Columns display the exact solution, approximate solution, log-scale loss history, and logarithmic relative $L^2$ error versus training iterations.
    }
    \label{fig:pacman}
\end{figure}

\end{example}


\section{Conclusions}

This work establishes a theoretical framework for employing low regularity neural networks for the numerical approximation of PDEs.
Emphasis is placed on hardware-oriented architectures, specifically those utilizing Heaviside activations and quantized weights. While our current implementations do not directly exploit the computational efficiency of binarized networks, this study provides a proof of concept for future deployment on specialized hardware. Numerical experiments demonstrate the efficacy of the method, even for problems characterized by low regularity solutions. Notably, binarized architectures maintain robust performance when utilizing large support test functions and surrogate gradient methods during optimization. This indicates that neither the test space saturation error nor the intrinsic error of the surrogate gradient approximation poses a significant limitation. Nevertheless, empirical results reveal a sensitivity to the hyperparameters $h$, $N_x$, $N_y$, and $N_b$. Deriving their optimal scaling remain as an open problem.

\subsection*{Declaration of AI-Assisted Technologies in Manuscript Preparation}
The authors acknowledge the use of AI tools to improve the English phrasing and grammar,  assuming complete responsibility for all the scientific content.

\bibliographystyle{abbrv}
\bibliography{main.bib}

@Inbook{Balinsky2015,
author="Balinsky, Alexander A.
and Evans, W. Desmond
and Lewis, Roger T.",
title="Hardy's Inequality on Domains",
bookTitle="The Analysis and Geometry of Hardy's Inequality",
year="2015",
publisher="Springer International Publishing",
address="Cham",
pages="77--134"
}

@article{kong2025expressivity,
  title={On the expressivity of deep Heaviside networks},
  author={Kong, Insung and Chen, Juntong and Langer, Sophie and Schmidt-Hieber, Johannes},
  journal={arXiv preprint arXiv:2505.00110},
  year={2025}
}

@inproceedings{umuroglu2017finn,
  title={Finn: A framework for fast, scalable binarized neural network inference},
  author={Umuroglu, Yaman and Fraser, Nicholas J and Gambardella, Giulio and Blott, Michaela and Leong, Philip and Jahre, Magnus and Vissers, Kees},
  booktitle={Proceedings of the 2017 ACM/SIGDA international symposium on field-programmable gate arrays},
  pages={65--74},
  year={2017}
}

@book{G,
 author = {Grisvard, P.},
 title = {Elliptic problems in nonsmooth domains},
 fseries = {Monographs and Studies in Mathematics},
 series = {Monogr. Stud. Math.},
 volume = {24},
 isbn = {0-273-08647-2},
 year = {1985},
 publisher = {Pitman, Boston, MA},
 language = {English},
 zbMATH = {4138299},
 Zbl = {0695.35060}
}

@book{LM,
 author = {Lions, J. L. and Magenes, E.},
 title = {Non-homogeneous boundary value problems and applications. {Vol}. {I,II,III}. {Translated} from the {French} by {P}. {Kenneth}},
 fseries = {Grundlehren der Mathematischen Wissenschaften},
 series = {Grundlehren Math. Wiss.},
 issn = {0072-7830},
 volume = {},
 year = {1972},
 publisher = {Springer, Cham},
 language = {English},
 zbMATH = {3353865},
 Zbl = {0223.35039}
}

@incollection{GPP,
 author = {Glowinski, R. and Pan, T.-W. and Periaux, J.},
 title = {Fictitious domain method for the {Dirichlet} problem and its generalization to some flow problems},
 booktitle = {Finite elements in fluids: new trends and applications. Proceedings of the 8th international conference, Barcelona, Spain, September 20-23, 1993. Part I},
 isbn = {84-87867-30-8; 84-87867-29-4},
 pages = {347--368},
 year = {1993},
 publisher = {Barcelona: Centro Internacional de M{\'e}todos Num{\'e}ricos en Ingenier{\'{\i}}a; Swansea: Pineridge Press},
 language = {English},
 zbMATH = {1023579},
 Zbl = {0876.76060}
}

@book{Evans,
 author = {Evans, Lawrence C.},
 title = {Partial differential equations},
 edition = {2nd ed.},
 fseries = {Graduate Studies in Mathematics},
 series = {Grad. Stud. Math.},
 issn = {1065-7339},
 volume = {19},
 isbn = {978-0-8218-4974-3; 978-1-4704-6942-9; 978-1-4704-1144-2},
 year = {2010},
 publisher = {Providence, RI: American Mathematical Society (AMS)},
 language = {English},
 zbMATH = {5681750},
 Zbl = {1194.35001}
}

@article{NikThom,
 author = {Nikolov, Nikolai and Thomas, Pascal J.},
 title = {Boundary regularity for the distance functions, and the eikonal equation},
 fjournal = {The Journal of Geometric Analysis},
 journal = {J. Geom. Anal.},
 issn = {1050-6926},
 volume = {35},
 number = {8},
 pages = {8},
 note = {Id/No 230},
 year = {2025},
 language = {English},
 doi = {10.1007/s12220-025-02064-7},
 zbMATH = {8072166},
 Zbl = {1570.35101}
}

@book{MaSch,
 author = {Maz'ya, Vladimir and Schmidt, Gunther},
 title = {Approximate approximations},
 fseries = {Mathematical Surveys and Monographs},
 series = {Math. Surv. Monogr.},
 issn = {0076-5376},
 volume = {141},
 isbn = {978-0-8218-4203-4},
 year = {2007},
 publisher = {Providence, RI: American Mathematical Society (AMS)},
 language = {English},
 zbMATH = {5189731},
 Zbl = {1120.41013}
}

@article{quant1,
  title={Binarized neural networks: Training deep neural networks with weights and activations constrained to+ 1 or-1},
  author={Courbariaux, Matthieu and Hubara, Itay and Soudry, Daniel and El-Yaniv, Ran and Bengio, Yoshua},
  journal={arXiv preprint arXiv:1602.02830},
  year={2016}
}

@article{hdim1,
  title={Tackling the curse of dimensionality with physics-informed neural networks},
  author={Hu, Zheyuan and Shukla, Khemraj and Karniadakis, George Em and Kawaguchi, Kenji},
  journal={Neural Networks},
  volume={176},
  pages={106369},
  year={2024},
  publisher={Elsevier}
}

@article{hdim2,
  title={Bias-variance trade-off in physics-informed neural networks with randomized smoothing for high-dimensional PDEs},
  author={Hu, Zheyuan and Yang, Zhouhao and Wang, Yezhen and Karniadakis, George E and Kawaguchi, Kenji},
  journal={SIAM Journal on Scientific Computing},
  volume={47},
  number={4},
  pages={C846--C872},
  year={2025},
  publisher={SIAM}
}

@article{hdim3,
  title={Hutchinson trace estimation for high-dimensional and high-order physics-informed neural networks},
  author={Hu, Zheyuan and Shi, Zekun and Karniadakis, George Em and Kawaguchi, Kenji},
  journal={Computer Methods in Applied Mechanics and Engineering},
  volume={424},
  pages={116883},
  year={2024},
  publisher={Elsevier}
}

@article{opschoor2024first,
  author  = {J. A. A. Opschoor and P. C. Petersen and C. Schwab},
  title   = {First order system least squares neural networks},
  journal = {arXiv preprint},
  volume  = {arXiv:2409.20264},
  year    = {2024}
}

@incollection{pinkus,
  author    = {A. Pinkus},
  title     = {Approximation theory of the {MLP} model in neural networks},
  booktitle = {Acta Numerica},
  volume    = {8},
  pages     = {143--195},
  publisher = {Cambridge University Press},
  year      = {1999}
}

@article{FNM,
  author  = {J. Xu},
  title   = {Finite neuron method and convergence analysis},
  journal = {Commun. Comput. Phys.},
  volume  = {28},
  pages   = {1707--1745},
  year    = {2020}
}

@article{PINN,
  author  = {M. Raissi and P. Perdikaris and G. E. Karniadakis},
  title   = {Physics-informed neural networks: A deep learning framework for solving forward and inverse problems involving nonlinear partial differential equations},
  journal = {J. Comput. Phys.},
  volume  = {378},
  pages   = {686--707},
  year    = {2019}
}

@article{liu2,
  author  = {M. Liu and Z. Cai},
  title   = {Adaptive two-layer {ReLU} neural network: II. {R}itz approximation to elliptic {PDE}s},
  journal = {Comput. Math. Appl.},
  volume  = {113},
  pages   = {103--116},
  year    = {2022}
}

@article{liu1,
  author  = {M. Liu and Z. Cai and J. Chen},
  title   = {Adaptive two-layer {ReLU} neural network: I. Best least-squares approximation},
  journal = {Comput. Math. Appl.},
  volume  = {113},
  pages   = {34--44},
  year    = {2022}
}

@inproceedings{weinan2022some,
  author    = {W. E and S. Wojtowytsch},
  title     = {Some observations on high-dimensional partial differential equations with {B}arron data},
  booktitle = {Proc. Math. Sci. Mach. Learn.},
  pages     = {253--269},
  year      = {2022}
}

@article{wojtowytsch2020can,
  author  = {S. Wojtowytsch and W. E},
  title   = {Can shallow neural networks beat the curse of dimensionality? A mean field training perspective},
  journal = {IEEE Trans. Artif. Intell.},
  volume  = {1},
  number  = {2},
  pages   = {121--129},
  year    = {2020}
}

@article{E18,
  author  = {W. E and B. Yu},
  title   = {The deep {R}itz method: A deep learning-based numerical algorithm for solving variational problems},
  journal = {Commun. Math. Stat.},
  volume  = {6},
  number  = {1},
  pages   = {1--12},
  year    = {2018}
}

@article{he20,
  author  = {C. He and X. Hu and L. Mu},
  title   = {A mesh-free method using piecewise deep neural network for elliptic interface problems},
  journal = {J. Comput. Appl. Math.},
  volume  = {412},
  pages   = {114358},
  year    = {2022}
}

@article{DGM18,
  author  = {J. Sirignano and K. Spiliopoulos},
  title   = {DGM: A deep learning algorithm for solving partial differential equations},
  journal = {J. Comput. Phys.},
  volume  = {375},
  pages   = {1339--1364},
  year    = {2018}
}

@article{DeepFOSLS,
  author  = {F. Bersetche and J. P. Borthagaray},
  title   = {A deep first-order system least squares method for solving elliptic {PDEs}},
  journal = {Comput. Math. Appl.},
  volume  = {129},
  pages   = {136--150},
  year    = {2023}
}
\end{document}